\documentclass[11pt]{amsart}
\usepackage{amsmath,amssymb,amsthm}
\usepackage[margin=1.15in]{geometry}
\usepackage{enumerate}
\usepackage{fancyvrb}
\usepackage{hyperref}

\newtheorem{theorem}{Theorem}[section]
\newtheorem{lemma}[theorem]{Lemma}
\newtheorem{proposition}[theorem]{Proposition}
\newtheorem{corollary}[theorem]{Corollary}
\newtheorem*{thmA}{Theorem A}
\newtheorem*{corB}{Corollary B}
\newtheorem*{thmC}{Theorem C}
\newtheorem*{thmD}{Theorem D}
\newtheorem*{thmE}{Theorem E}
\theoremstyle{definition}

\newtheorem{remark}[theorem]{Remark}
\newtheorem{hypothesis}[theorem]{Hypothesis}

\newcommand{\F}{\mathbb{F}}
\newcommand{\GF}{\mathrm{GF}}
\newcommand{\X}{\widetilde{\mathcal{X}}_p}
\newcommand{\Zh}{\widehat{Z}}
\newcommand{\Oh}{\widehat{\Omega}}
\newcommand{\Aut}{\mathrm{Aut}}
\newcommand{\Out}{\mathrm{Out}}
\newcommand{\Inn}{\mathrm{Inn}}
\newcommand{\Sym}{\mathrm{Sym}}
\newcommand{\Alt}{\mathrm{Alt}}
\newcommand{\Syl}{\mathrm{Syl}}
\newcommand{\Sp}{\mathrm{Sp}}
\newcommand{\GL}{\mathrm{GL}}
\newcommand{\SL}{\mathrm{SL}}
\newcommand{\PSL}{\mathrm{PSL}}
\newcommand{\PSU}{\mathrm{PSU}}
\newcommand{\PSp}{\mathrm{PSp}}
\newcommand{\Sz}{\mathrm{Sz}}
\newcommand{\End}{\mathrm{End}}
\newcommand{\Gal}{\mathrm{Gal}}
\newcommand{\Irr}{\mathrm{Irr}}
\newcommand{\cl}{\mathrm{cl}}

\newcommand{\Fs}{F^{*}}

\begin{document}

\title[Large orbits of $p$-groups on conjugacy classes]{Large orbits of $p$-groups on conjugacy classes of arbitrary finite groups}

\author{Thomas Michael Keller}
\address{Department of Mathematics, Texas State University, 601 University Drive, San Marcos, TX 78666, United States of America}
\email{keller@txstate.edu}


\begin{abstract}
Let $P$ be a $p$-group that acts faithfully and coprimely on a finite group $H$. In a previous paper it was shown that if $H$ is solvable, then there exists a conjugacy class $c$ of $H$ whose stabilizer in $P$ has order at most $|P|^{1/2}$, and that this bound is best possible. In this paper we remove the hypothesis that $H$ is solvable. As a consequence, there always exists an irreducible character of $H$ whose inertia group in $P$ has order at most $|P|^{1/2}$. The proof depends on the classification of finite simple groups. Its main new ingredient is a result on quasiprimitive linear actions of $p$-nilpotent groups which shows that the elements acting as field automorphisms can be ignored when one looks for vectors with small stabilizers. We also present a new family of examples, with nonsolvable $H$, in which the bound is attained.
\end{abstract}

\maketitle

\section{Introduction}

One important aspect of the study of finite group actions is the existence of large orbits, or equivalently, of points with small stabilizers. A classical result of this type is due to D.~S.~Passman \cite{Pa} and I.~M.~Isaacs \cite{Is99}: if $P$ is a nontrivial $p$-group that acts faithfully on a group $H$ of order not divisible by $p$, then there exists an element $x\in H$ such that $|C_P(x)|\le (|P|/p)^{1/p}$.

In \cite{K} the corresponding question was studied for the action of $P$ on the set $\cl(H)$ of conjugacy classes of $H$. If $c\in\cl(H)$, we write $C_P(c)$ for the stabilizer of $c$ in $P$. It was shown in \cite{K} that for classes the bound of Passman and Isaacs fails for every odd prime $p$, and that the correct bound is $|P|^{1/2}$, at least when $H$ is solvable. The purpose of this paper is to show that the solvability hypothesis is not needed.

\begin{thmA}
Let $P$ be a $p$-group that acts faithfully and coprimely on a finite group $H$. Then there exists a conjugacy class $c$ of $H$ such that
\[ |C_P(c)|\le |P|^{1/2}. \]
\end{thmA}

By \cite[Theorem B]{K}, this bound is best possible for every prime $p$. Since the coprime actions of $P$ on the classes and on the irreducible characters of $H$ are permutation isomorphic, we also obtain the following.

\begin{corB}
Let $P$ be a $p$-group that acts faithfully and coprimely on a finite group $H$. Then there exists $\chi\in\Irr(H)$ such that $|I_P(\chi)|\le |P|^{1/2}$, where $I_P(\chi)$ denotes the stabilizer of $\chi$ in $P$.
\end{corB}

Note that if $p=2$, then $|H|$ is odd, so that $H$ is solvable by the Feit--Thompson theorem and Theorem A is contained in \cite{K}. Thus only odd primes have to be considered. For odd $p$ our proof depends on the classification of finite simple groups, through the list of the simple groups whose order is not divisible by $p$ and the structure of their automorphism groups.

As in \cite{K}, Theorem A is obtained from a result on linear actions of $p$-nilpotent groups. We denote by $\X$ the class of all pairs $(G,V)$, where $G$ is a finite $p$-nilpotent group with $O_p(G)>1$ and $V$ is a finite faithful irreducible $G$-module over a field of prime order. This is the class $\mathcal X_p$ of \cite{K} without the hypothesis that $G$ is solvable. For a group $N$ acting on a set $V$ we write $\mu_p(N,V)=\min_{v\in V}|C_N(v)|_p$, and we refer to Section 2 for the notions of tight and rigid actions.

\begin{thmC}
Let $p$ be a prime, let $(G,V)\in\X$ and let $N$ be a normal subgroup of $G$. Then the following hold.
\begin{enumerate}[\rm (a)]
\item There exists $v\in V$ such that $|C_N(v)|_p\le |N|_p^{1/2}$.
\item Suppose that $|N|_p=p^{2a}$ for some $a\ge 1$ and that $|C_N(v)|_p\ge p^a$ for all $v\in V$. Then the vectors $v\in V$ with $|C_N(v)|_p=p^a$ lie in at least two different $N$-orbits.
\end{enumerate}
\end{thmC}

We now describe the new ideas. The reduction of Theorem A to Theorem C in \cite{K} uses the Fitting series of $H$. Here we replace the Fitting subgroup by the generalized Fitting subgroup. This leads to an additional problem for groups acting on a direct product $E$ of nonabelian simple groups, which is solved in Section 5: if $X$ is $p$-nilpotent with $\Inn(E)\le X\le\Aut(E)$, then there is an element $e\in E$ such that $C_X(e)$ is a $p'$-group. The proof needs a result on colourings of permutation groups (Lemma \ref{colouring}), for which two colours do not suffice.

The proof of Theorem C is by induction on $|G|$. If $V$ is not quasiprimitive, the arguments of \cite[Section 5]{K} apply without change. In the quasiprimitive case the proof in \cite{K} relies on the classification of Holt and Yang \cite{HY} and on computer calculations; there the goal is, in most cases, to find a vector $v$ such that $C_G(v)$ is a $p'$-group. For nonsolvable groups this strategy cannot work in general, as the following result shows.

\begin{thmE}
Let $p$ be a prime with $p\equiv 1\pmod 3$. Then there exist pairs $(G,V)\in\X$ such that $O_{p'}(G)$ is nonsolvable, $|C_G(v)|_p=p$ for every nonzero $v\in V$, and such that $G$ has a normal subgroup $N$ with $|N|_p=p^2$ and $|C_N(v)|_p=p=|N|_p^{1/2}$ for every nonzero $v\in V$.
\end{thmE}

In these examples, which are described in Section 12, $O_{p'}(G)$ is a central product of $\SL_2(q)$ and a cyclic group, and the $p$-elements outside $O_p(G)$ act on $V$ as field automorphisms. So the inequality in Theorem C(a) is attained also for nonsolvable groups. To deal with such situations we prove the following result, which is the main new tool of this paper. It is stated here in an informal way; see Theorem \ref{galois} for the precise statement.

\begin{thmD}
Let $(G,V)\in\X$, and let $\F$ be a subfield of $\End(V)$ which is normalized by $G$ and such that $p$ divides $|\F|-1$. Let $B=C_G(\F)$. If there exists a nonzero vector in $V$ which is not an eigenvector of any $p$-element of $B$ that is not a scalar, then Theorem C holds for $(G,V)$ and every normal subgroup $N$ of $G$.
\end{thmD}

Thus the elements of $G$ which act on $\F$ as nontrivial field automorphisms play no role. The proof of Theorem D uses the scalars of $p$-power order in $\F$, which need not belong to $G$, but which move the stabilizers of vectors in a controlled way. This idea is already implicit in the treatment of semilinear groups in \cite[Lemma 6.4]{K}. Theorem D reduces the quasiprimitive case to a counting problem for $\F$-linear elements, which we solve in Sections 8 to 10. For $p\ge 5$ this argument does not need that $G$ is nonsolvable, and so it gives a proof of the quasiprimitive case of \cite[Theorem C]{K} which does not use \cite{HY}.

The paper is organized as follows. In Section 2 we collect notation and auxiliary results. In Section 3 we reduce Theorem A to Theorem C and to the result on products of simple groups mentioned above, which is proved in Section 5. Section 4 deals with the imprimitive case of Theorem C. In Section 6 we describe the structure of the groups in the quasiprimitive case, and in Section 7 we prove Theorem D. Sections 8, 9 and 10 contain the counting arguments, and in Section 11 we complete the proofs of Theorems A and C and Corollary B. Section 12 contains the proof of Theorem E, and we conclude with some remarks in Section 13. Two small computer calculations are described in the appendix.

\section{Notation and auxiliary results}

Throughout, all groups and modules are finite, and we use the notation of \cite{K}. If $n$ is a positive integer and $p$ is a prime, then $n_p$ is the $p$-part of $n$, and $|X|_p=(|X|)_p$ for a group $X$. If a group $X$ acts on a set $\Omega$ and $\omega\in\Omega$, then $C_X(\omega)$ denotes the stabilizer of $\omega$ in $X$. A module $V$ for a group $G$ is called quasiprimitive if $V_N$ is homogeneous for all $N\unlhd G$. By $F(X)$, $E(X)$ and $\Fs(X)=F(X)E(X)$ we denote the Fitting subgroup, the layer and the generalized Fitting subgroup of a group $X$; recall that $C_X(\Fs(X))\le \Fs(X)$, and that the components of $X$ are its subnormal quasisimple subgroups.

Let $N$ be a group acting on a finite set $V$, and let $p$ be a prime. As in \cite[Definition 2.1]{K} we put $\mu_p(N,V)=\min_{v\in V}|C_N(v)|_p$, and we call $v\in V$ a \emph{minimizer} if $|C_N(v)|_p=\mu_p(N,V)$. The action of $N$ on $V$ is \emph{tight} if $|N|_p=p^{2a}$ for some $a\ge 1$ and $\mu_p(N,V)=p^a$, and it is \emph{rigid} if it is tight and the minimizers form a single $N$-orbit. With this notation, Theorem C states that for $(G,V)\in\X$ and $N\unlhd G$ we have $\mu_p(N,V)\le |N|_p^{1/2}$, and that the action of $N$ on $V$ is not rigid. In particular, Theorem C holds for $(G,V)$ if there is a vector $v\in V$ such that $C_G(v)$ is a $p'$-group.

We recall four results from \cite{K}.

\begin{lemma}[{\cite[Lemma 2.2]{K}}]\label{classes}
Let $P$ be a $p$-group acting on a $p'$-group $H$, and let $G=HP$ be the semidirect product. Let $x\in H$ and let $c=x^H$. Then $|C_P(c)|=|C_G(x)|_p$.
\end{lemma}

\begin{lemma}[{\cite[Lemma 2.3]{K}}]\label{product}
Let a group $N$ act on a vector space $V=V_1\oplus\dots\oplus V_m$, where the $V_j$ are $N$-invariant. Put $M_0=N$ and $M_j=C_N(V_1\oplus\dots\oplus V_j)$ for $j=1,\dots,m$. Let $v_j\in V_j$ for all $j$ and $v=v_1+\dots+v_m$. Then
\[ |C_N(v)|_p\le |M_m|_p\prod_{j=1}^m |C_{M_{j-1}/M_j}(v_j)|_p. \]
The same holds, with the same proof, if $V=V_1\times\dots\times V_m$ is a direct product of groups on which $N$ acts by automorphisms and the $V_j$ are $N$-invariant.
\end{lemma}

\begin{lemma}[{\cite[Lemma 2.5]{K}}]\label{semilinear}
Let $K=\GF(Q)$, let $\sigma$ be an automorphism of $K$ of order $s$, and let $y$ be a bijective $\sigma$-semilinear map of a $K$-vector space $W$ of dimension $m$. Then $|C_W(y)|\le Q^{m/s}$.
\end{lemma}

\begin{lemma}[{\cite[Lemma 2.6]{K}}]\label{arith}
Let $p$ be a prime, $Q_0$ a prime power and $j\ge 0$.
\begin{enumerate}[\rm (a)]
\item If $p$ divides $Q_0^{p^j}-1$, then $p$ divides $Q_0-1$.
\item If $p$ divides $Q_0-1$, then $p^j$ divides $(Q_0^{p^j}-1)/(Q_0-1)$.
\end{enumerate}
\end{lemma}

The next two lemmas are standard facts on coprime actions.

\begin{lemma}\label{fstar}
Let $\alpha$ be an automorphism of a group $X$ whose order is coprime to $|X|$. If $\alpha$ centralizes $\Fs(X)$, then $\alpha=1$.
\end{lemma}

\begin{proof}
We have $[\Fs(X),\alpha,X]=1$ and $[X,\Fs(X),\alpha]\le[\Fs(X),\alpha]=1$. By the three subgroups lemma, $[\alpha,X,\Fs(X)]=1$, so that $[X,\alpha]\le C_X(\Fs(X))\le \Fs(X)$. Hence $[X,\alpha,\alpha]=1$, and since the action is coprime, $[X,\alpha]=[X,\alpha,\alpha]=1$.
\end{proof}

\begin{lemma}\label{complements}
Let $L$ be a $p'$-group and let $\alpha$ be an automorphism of $L$ of order $p$. Then the following hold in the semidirect product $L\langle\alpha\rangle$.
\begin{enumerate}[\rm (a)]
\item Every element of order $p$ in the coset $\alpha L$ is conjugate to $\alpha$ by an element of $L$.
\item For every prime $s$, there is an $\alpha$-invariant Sylow $s$-subgroup of $L$.
\item An element $x\in L$ can be written as $x=\alpha(g)g^{-1}$ with $g\in L$ if and only if $\alpha^{-1}x$ has order $p$. In particular, this is the case if $x$ lies in an $\alpha$-invariant abelian subgroup $T$ of $L$ with $C_T(\alpha)=1$.
\end{enumerate}
Here $\alpha(g)=\alpha g\alpha^{-1}$.
\end{lemma}

\begin{proof}
(a) If $y\in\alpha L$ has order $p$, then $\langle y\rangle$ is a complement of $L$ in $L\langle\alpha\rangle$, and so $\langle y\rangle=\langle\alpha\rangle^g$ for some $g\in L$ by the Schur--Zassenhaus theorem. As $y\in\alpha L$, we get $y=\alpha^g$.

(b) Let $Q\in\Syl_s(L)$. By the Frattini argument, $N_{L\langle\alpha\rangle}(Q)$ contains an element of $\alpha L$, and hence, taking a suitable power, an element of order $p$ in $\alpha L$. By (a) this element is $\alpha^g$ for some $g\in L$, and so $\alpha$ normalizes $Q^{g^{-1}}$.

(c) We have $x=\alpha(g)g^{-1}$ if and only if $\alpha^{-1}x=g\alpha^{-1}g^{-1}$, so the first statement follows from (a), applied to $\alpha^{-1}$. If $x\in T$ and $\alpha^{-1}$ acts on $T$, then $(\alpha^{-1}x)^p=\prod_{i=1}^{p}\alpha^{-i}(x)$, and this product is an element of $C_T(\alpha)=1$.
\end{proof}

We will also use Zsigmondy's theorem: if $\ell$ is a prime and $n\ge 3$ with $(\ell,n)\ne(2,6)$, then there is a prime $r$ such that $\ell$ has order $n$ modulo $r$. We call such a prime a \emph{Zsigmondy prime} of $\ell^n-1$. Note that then $n$ divides $r-1$.

\section{The reduction of Theorem A}

In this section we reduce Theorem A to Theorem C(a) and to the following result, which will be proved in Section 5.

\begin{theorem}\label{layer}
Let $p$ be an odd prime, and let $E=S_1\times\dots\times S_m$ be a direct product of nonabelian simple groups whose orders are not divisible by $p$. Let $X$ be a $p$-nilpotent group with $\Inn(E)\le X\le\Aut(E)$. Then there exists $e\in E$ such that $C_X(e)$ is a $p'$-group.
\end{theorem}

The following lemma replaces \cite[Lemma 4.1]{K}.

\begin{lemma}\label{reduction}
Let $p$ be an odd prime, and assume that $\mu_p(G,V)\le|G|_p^{1/2}$ for all $(G,V)\in\X$. Let $G$ be a $p$-nilpotent group such that a Sylow $p$-subgroup $P$ of $G$ acts faithfully on $H=O_{p'}(G)$. Then there exists $x\in H$ such that $|C_G(x)|_p\le|G|_p^{1/2}$.
\end{lemma}

\begin{proof}
We argue by induction on $|G|$. If $P=1$, there is nothing to prove, so let $P>1$.

\textbf{Step 1.} $O_p(G)=1$ and $\Fs(G)\le H$.

As in \cite{K}, $O_p(G)$ centralizes $H$, and so $O_p(G)=1$ by hypothesis. Hence $F(G)\le H$. If $L$ is a component of $G$, then $LH/H$ is a perfect subnormal subgroup of the $p$-group $G/H$, so that $L\le H$. Thus $\Fs(G)\le H$.

\textbf{Step 2.} We may assume that $\Phi(G)=1$.

Let $\Phi=\Phi(G)$ and suppose that $\Phi>1$. By Step 1, $\Phi\le F(G)\le H$. We show that $G/\Phi$ satisfies the hypotheses of the lemma. Clearly $G/\Phi$ is $p$-nilpotent with $O_{p'}(G/\Phi)=H/\Phi$. Let $R/\Phi=O_p(G/\Phi)$ and $S\in\Syl_p(R)$. By the Frattini argument, $G=RN_G(S)=\Phi N_G(S)$, and thus $G=N_G(S)$. So $S\le O_p(G)=1$ and $O_p(G/\Phi)=1$. By the argument of Step 1, $\Fs(G/\Phi)\le H/\Phi$. Now if $y\in P$ acts trivially on $H/\Phi$, then $y\Phi$ centralizes $\Fs(G/\Phi)$, so that $y\Phi\in\Fs(G/\Phi)\le H/\Phi$ and thus $y\in H\cap P=1$. Hence by induction there is $x\in H$ with $|C_{G/\Phi}(x\Phi)|_p\le|G/\Phi|_p^{1/2}=|G|_p^{1/2}$. Since $C_G(x)\Phi/\Phi\le C_{G/\Phi}(x\Phi)$ and $\Phi$ is a $p'$-group, we are done.

\textbf{Step 3.} The structure of $\Fs(G)$.

Put $F=F(G)$ and $E=E(G)$. Since $\Phi(F)\le\Phi(G)=1$, the group $F$ is abelian. Since $E$ is perfect, $Z(E)\le\Phi(E)\le\Phi(G)=1$. Hence $\Fs(G)=F\times E$, and $E$ is a direct product of nonabelian simple groups, whose orders are not divisible by $p$ by Step 1.

\textbf{Step 4.} The case $E\ne 1$.

Let $M=C_G(E)$. Then $M$ is a normal $p$-nilpotent subgroup of $G$ with $O_{p'}(M)=C_H(E)$, and $M<G$. If $y\in M\cap P$ centralizes $C_H(E)$, then $y$ centralizes $F$ and $E$, so that $y\in\Fs(G)\le H$ and $y=1$. Hence by induction there exists $f\in C_H(E)$ with $|C_M(f)|_p\le|M|_p^{1/2}$. Since $C_H(E)\cap E=Z(E)=1$ and both subgroups are normal in $G$, we have $C_G(fe)=C_G(f)\cap C_G(e)$ for every $e\in E$. Now $G/M$ is a $p$-nilpotent group which acts faithfully on $E$ and contains $EM/M\cong\Inn(E)$. By Theorem \ref{layer} there exists $e\in E$ such that $C_G(e)M/M$ is a $p'$-group. For $x=fe\in H$ we obtain
\[ |C_G(x)|_p=|C_G(x)M/M|_p\,|C_G(x)\cap M|_p\le|C_M(f)|_p\le|M|_p^{1/2}\le|G|_p^{1/2}. \]

\textbf{Step 5.} The case $E=1$.

Now $\Fs(G)=F$ and $C_G(F)=F$. By a theorem of Gasch\"utz (see \cite[III, Hilfssatz 4.4 and Satz 4.5]{Hu}), which is proved there for arbitrary finite groups, $F$ is a direct product of minimal normal subgroups of $G$, and, being abelian with $F\cap\Phi(G)=1$, it has a complement $L$ in $G$. Thus $L\cong G/F$ acts faithfully on $F$, and $F\ne 1$.

Suppose first that $O_p(L)=1$. A Sylow $p$-subgroup $P_L$ of $L$ which centralizes $O_{p'}(L)$ is a direct factor of $C_L(O_{p'}(L))$, hence contained in $O_p(L)$. So $P_L$ acts faithfully on $O_{p'}(L)\le H$, and by induction there exists $y\in O_{p'}(L)$ with $|C_L(y)|_p\le|L|_p^{1/2}=|G|_p^{1/2}$. Since $C_G(y)F/F$ is contained in the centralizer of $yF$ in $G/F\cong L$, we get $|C_G(y)|_p\le|C_L(y)|_p$, and we are done.

Now suppose that $O_p(L)>1$. Then Step 5 of the proof of \cite[Lemma 4.1]{K} can be repeated word for word; it only uses that $F$ is an abelian direct product of minimal normal subgroups of $G$ with $C_G(F)=F$, and the inequality $\mu_p(\bar L,V_1)\le|\bar L|_p^{1/2}$ for a pair $(\bar L,V_1)$ which now lies in $\X$.
\end{proof}

\section{The imprimitive case}

In this section and in Sections 6 to 10 we prove Theorem C for odd primes $p$ by induction on $|G|$. As in \cite{K} we use the following inductive hypothesis.

\begin{hypothesis}\label{ind}
Theorem C holds for all pairs $(\tilde G,\tilde V)\in\X$ with $|\tilde G|<|G|$ and all normal subgroups of $\tilde G$.
\end{hypothesis}

\begin{proposition}\label{imprimitive}
Let $(G,V)\in\X$ such that $V$ is not quasiprimitive, and assume Hypothesis \ref{ind}. Then Theorem C holds for $(G,V)$.
\end{proposition}

\begin{proof}
This is proved in \cite[Section 5]{K} for solvable $G$, and we only have to check that solvability is not used there. Put $K=O_{p'}(G)$ and $Z=O_p(G)$. The proof of \cite[Lemma 5.3]{K} uses that $G/D$ is a $p$-group for a normal subgroup $D$ containing $K$, and \cite[Lemma 2.4]{K}, which is a statement on $p$-groups. The proof of \cite[Proposition 5.4]{K} uses Clifford's theorem, the fact that $[K,Z]=1$, and that $KQ\unlhd G$ whenever $Q$ is a normal subgroup of a Sylow $p$-subgroup of $G$. The proof of \cite[Theorem 5.5]{K} uses Sylow's theorem in a normal subgroup of the stabilizer of a system of imprimitivity, and the proofs of \cite[Lemma 5.6 and Theorem 5.7]{K} are counting arguments based on Lemma \ref{product}. Finally, the pairs $(J/C_J(U_1),U_1)$ and $(I/C_I(V_1),V_1)$ to which the inductive hypothesis is applied in \cite{K} belong to $\X$, because subgroups and factor groups of $p$-nilpotent groups are $p$-nilpotent.
\end{proof}

\section{Groups acting on products of simple groups}

In this section we prove Theorem \ref{layer}. Throughout, $p$ is an odd prime. We need some facts on simple groups whose order is not divisible by $p$, which depend on the classification of finite simple groups.

\begin{lemma}\label{out}
Let $S$ be a nonabelian simple group with $p\nmid|S|$. Then the Sylow $p$-subgroups of $\Out(S)$ are cyclic. If $p$ divides $|\Out(S)|$, then $S$ is a group of Lie type over a field $\GF(q)$ with $q=q_0^p$, an element of order $p$ in $\Out(S)$ is the image of a field automorphism $\varphi$ of order $p$, and $C_S(\varphi)$ is isomorphic to a subgroup of the group $G_{ad}(q_0)$ of inner-diagonal automorphisms of the group of the same type as $S$ over $\GF(q_0)$.
\end{lemma}

\begin{proof}
If $S$ is alternating or sporadic, then $\Out(S)$ is a $2$-group. So let $S$ be of Lie type over $\GF(q)$. By \cite[Theorem 2.5.12(b)]{GLS}, $\Out(S)$ is a split extension of $\mathrm{Outdiag}(S)$ by $\Phi_S\Gamma_S$, where $\Phi_S$ is the cyclic group of field automorphisms and $\Gamma_S$ the group of graph automorphisms. By part (c) of that theorem, every prime divisor of $|\mathrm{Outdiag}(S)|$ divides $|S|$, so $p\nmid|\mathrm{Outdiag}(S)|$. By parts (d)--(f), $\Phi_S\Gamma_S$ is either cyclic, or a direct product $\Phi_S\times\Gamma_S$ with $\Gamma_S$ a group of symmetries of the Dynkin diagram, of order dividing $2$, or $6$ for type $D_4$; in the latter case $3$ divides $|S|$. Thus a Sylow $p$-subgroup of $\Out(S)$ is isomorphic to one of the cyclic group $\Phi_S$. The last statement is \cite[Proposition 4.9.1]{GLS}.
\end{proof}

We write the order of a group $G_{ad}(q)$ of Lie type as $q^N\prod_k\Phi_k(q)^{a_k}$, where $\Phi_k$ is the $k$-th cyclotomic polynomial, and we let $I(G)$ be the set of all $k$ with $a_k>0$; for the Suzuki and Ree groups we use the order formula as a polynomial in $q$, so that for instance $I=\{1,4\}$ for $\Sz(q)$. The sets $I(G)$ are listed in Table \ref{indices}.

\begin{table}[h]
\begin{tabular}{l|l}
Type & $I(G)$\\ \hline
$A_n$ & $1,\dots,n+1$\\
${}^2A_n$ & all $k$ with $k^*\le n+1$, where $k^*=2k$, $k/2$, $k$ for $k$ odd, $k\equiv 2\ (4)$, $4\mid k$\\
$B_n$, $C_n$ & all $k\le n$ and all even $k\le 2n$\\
$D_n$ & all divisors of $2i$ for $i\le n-1$, and of $n$\\
${}^2D_n$ & all divisors of $2i$ for $i\le n-1$, and the divisors of $2n$ not dividing $n$\\
$G_2$; ${}^3D_4$ & $1,2,3,6$; $1,2,3,6,12$\\
$F_4$ & $1,2,3,4,6,8,12$\\
$E_6$; ${}^2E_6$ & $1,2,3,4,5,6,8,9,12$; $1,2,3,4,6,8,10,12,18$\\
$E_7$ & $1,2,3,4,5,6,7,8,9,10,12,14,18$\\
$E_8$ & $1,\dots,10,12,14,15,18,20,24,30$\\
${}^2B_2$; ${}^2G_2$; ${}^2F_4$ & $1,4$; $1,2,6$; $1,2,4,6,12$
\end{tabular}
\medskip
\caption{The cyclotomic indices of the groups of Lie type.}\label{indices}
\end{table}

\begin{lemma}\label{zsig}
Let $p\ge 5$ and let $S$ be a simple group of Lie type over $\GF(q)$, where $q=\ell^f=q_0^p$ for a prime $\ell$, and suppose that $p\nmid|S|$. Let $k\in I(G)$ and let $r_k$ be a Zsigmondy prime of $\ell^{fk}-1$. Then $r_k$ divides $|S|$, but $r_k$ does not divide $|G_{ad}(q_0)|$. Moreover, $p$ divides $r_k-1$.
\end{lemma}

\begin{proof}
Note that $fk\ge 5$ and $fk\ne 6$, as $p$ divides $f$, so that $r_k$ exists, and $fk$ divides $r_k-1$. The order of $q$ modulo $r_k$ is $k$, and so $r_k$ divides $\Phi_k(q)$ and hence $|G_{ad}(q)|$; since $r_k>fk\ge 5$, it is easy to see that $r_k$ also divides $|S|$.

We first show that no element of $I(G)$ is divisible by $p$. Inspection of Table \ref{indices} shows that if $I(G)$ contains a multiple of $p$, then $I(G)$ contains all divisors of $p-1$. For example, for type $B_n$ a multiple of $p$ in $I(G)$ forces $p\le n$, and all $k\le n$ belong to $I(G)$; for type $E_8$ only $p\in\{5,7\}$ occur, and $1,2,3,4,6\in I(G)$. But then $I(G)$ contains the order $o$ of $q$ modulo $p$, so that $p$ divides $\Phi_o(q)$ and thus $|G_{ad}(q)|$, and then $p$ divides $|S|$ as in the proof of Lemma \ref{out}, a contradiction.

Now the order of $q_0=\ell^{f/p}$ modulo $r_k$ is $fk/\gcd(fk,f/p)=pk$. If $r_k$ divides $|G_{ad}(q_0)|$, then $r_k$ divides $\Phi_j(q_0)$ for some $j\in I(G)$, and then $j=pk\,r_k^{\,a}$ for some $a\ge 0$. This contradicts the previous paragraph.
\end{proof}

Now let $E$ and $X$ be as in Theorem \ref{layer}. For $i=1,\dots,m$ let $Y_i$ be the group of automorphisms of $S_i$ induced by $N_X(S_i)$, so that $\Inn(S_i)\le Y_i\le\Aut(S_i)$. We call an element of $S_i$ \emph{good} if it is not fixed by any element of order $p$ of $Y_i$. The $X$-orbits of good elements are called \emph{types}; note that the number of types contained in $\bigcup_j S_j^{\,x}$, $x\in X$, equals the number of orbits of $Y_i$ on the good elements of $S_i$.

\begin{lemma}\label{onecomponent}
Let $S$ be a nonabelian simple group with $p\nmid|S|$, and let $\Inn(S)\le Y\le\Aut(S)$. Then $S$ has good elements with respect to $Y$. More precisely, if $p$ divides $|Y|$, then there is a subgroup $Y_0$ of order $p$ of $Y$ such that the elements of $S$ which are not good are the elements of $\bigcup_{s\in S}C_S(Y_0)^s$.
\end{lemma}

\begin{proof}
If $p\nmid|Y|$, all elements are good. Otherwise, by Lemma \ref{out}, $Y/S$ has a unique subgroup of order $p$. Let $Y_0\le Y$ be of order $p$. If $Y_1$ is any subgroup of order $p$ of $Y$, then $Y_1S=Y_0S$, and by Sylow's theorem $Y_1$ is conjugate to $Y_0$ by an element of $S$. This proves the second statement, and the first one follows since $C_S(Y_0)<S$ and a group is not the union of the conjugates of a proper subgroup.
\end{proof}

\begin{lemma}\label{threetypes}
Let $S$ and $Y$ be as in Lemma \ref{onecomponent}.
\begin{enumerate}[\rm (a)]
\item If $p\ge 5$, then $Y$ has at least three orbits on the good elements of $S$.
\item If $p=3$, then $Y$ has at least three orbits on the good elements of $S$, unless $S\cong\Sz(8)$, $Y=\Aut(S)$, and the six conjugacy classes of elements of order $7$ and $13$ of $S$ are permuted by $Y$ in two orbits of length $3$. In the exceptional case these six classes consist of the good elements.
\end{enumerate}
\end{lemma}

\begin{proof}
If $p\nmid|Y|$, all elements of $S$ are good, and $S$ has elements of at least four different orders. So let $p$ divide $|Y|$. By Lemmas \ref{out} and \ref{onecomponent}, $S$ is of Lie type over $\GF(q)$ with $q=\ell^f=q_0^p$, and every element of $S$ whose order is not the order of an element of $G_{ad}(q_0)$ is good. Since elements of different orders lie in different $Y$-orbits, it suffices to find three such orders.

(a) By Lemma \ref{zsig}, each $k\in I(G)$ gives a prime $r_k$ which is the order of a good element. This yields three orders unless $|I(G)|=2$, that is, unless $S$ is $\PSL_2(q)$ or $\Sz(q)$. If $S=\Sz(q)$, then $S$ has cyclic subgroups of the orders $q-1$ and $q\pm\sqrt{2q}+1$, and each of these three numbers is larger than the order of any element of $\Sz(q_0)$. If $S=\PSL_2(q)$ and $d=\gcd(2,q-1)$, then $S$ has cyclic subgroups of the orders $(q-1)/d$ and $(q+1)/d$, which exceed all element orders of $\mathrm{PGL}_2(q_0)$. So we obtain a third order unless $(q-1)/d=r_1$ and $(q+1)/d=r_2$ are both primes. This is impossible: if $q$ is even, then $3$ divides $q+1>3$ because $f$ is odd; and if $q$ is odd, then $(q-1)/2$ is divisible by $(\ell-1)/2$, which forces $\ell=3$, and then $(q+1)/2$ is even and larger than $2$.

(b) Let $p=3$. Then $S=\Sz(q)$ with $q=2^f$ and $3\mid f$. If $q\ge 2^9$, let $r_1$ and $r_4$ be Zsigmondy primes of $2^f-1$ and $2^{4f}-1$. Neither divides $|\Sz(q_0)|$, and $q-1\ne r_1$ since $q_0-1>1$ divides $q-1$. So elements of the orders $r_1$, $r_4$ and $q-1$ are good. Finally let $q=8$, so that $C_S(Y_0)\cong\Sz(2)$ is a Frobenius group of order $20$ and $Y=\Aut(S)$. The group $\Sz(8)$ has eleven conjugacy classes, with representatives of the orders $1,2,4,4,5,7,7,7,13,13,13$, and $\Sz(2)$ meets exactly the classes of elements of the orders $1,2,4,5$. Hence the good elements are those of order $7$ or $13$, and the assertion follows.
\end{proof}

The next lemma is the combinatorial heart of this section.

\begin{lemma}\label{colouring}
Let $\bar X\le\Sym(\Omega)$ be a transitive $p$-nilpotent permutation group, where $p$ is odd. Then there exists a colouring $c:\Omega\to\{1,2,3\}$ with the following properties.
\begin{enumerate}[\rm (i)]
\item No element of order $p$ of $\bar X$ preserves $c$.
\item The three colour classes $c^{-1}(1)$, $c^{-1}(2)$, $c^{-1}(3)$ do not all have the same size.
\end{enumerate}
\end{lemma}

We remark that two colours do not suffice. Indeed, let $p=3$ and $\bar X=A_4\wr C_3$, acting on $12$ points in three blocks of size four. A $2$-colouring in which some block is not split into two classes of size two is preserved by a $3$-cycle of the corresponding factor $A_4$. If every block is split $2+2$, and $S_j$ is the set of points of colour $1$ in the $j$-th block, choose $g_1,g_2\in A_4$ mapping $S_1$ to $S_2$ and $S_2$ to $S_3$; then the element $(g_1,g_2,(g_1g_2)^{-1})t$, where $t$ permutes the blocks cyclically, has order $3$ and preserves the colouring.

\begin{proof}[Proof of Lemma \ref{colouring}]
We call a colouring with properties (i) and (ii) good, and we argue by induction on $n=|\Omega|$.

\textbf{Step 1.} The case that $\bar X$ is not primitive.

Let $\Sigma=\{B_1,\dots,B_t\}$ be a nontrivial system of blocks. The group induced by the stabilizer of $B_1$ on $B_1$ and the group $\bar X^\Sigma$ induced on $\Sigma$ are transitive and $p$-nilpotent, so by induction there are good colourings $c_1$ of $B_1$ and $\gamma$ of $\Sigma$. Let $\mathbf n=(n_1,n_2,n_3)$ and $\mathbf u=(u_1,u_2,u_3)$ be the vectors of the sizes of their colour classes. Let $\tau$ be the cyclic permutation $(1\,2\,3)$ of the colours. Since $\mathbf n$ is not constant, the three colourings $\tau^i\circ c_1$, $i=0,1,2$, have pairwise different vectors of class sizes. For each $j$ fix $g_j\in\bar X$ with $B_1^{g_j}=B_j$, and define $c$ on $B_j$ by $c(\omega^{g_j})=\tau^{\gamma(B_j)-1}(c_1(\omega))$ for $\omega\in B_1$.

Let $y\in\bar X$ be of order $p$ and suppose that $y$ preserves $c$. If $y$ acts nontrivially on $\Sigma$, then, since $y$ maps blocks to blocks with the same vector of class sizes, $y$ preserves $\gamma$, and it induces an element of order $p$ on $\Sigma$, a contradiction. If $y$ fixes every block, then $y$ acts nontrivially on some $B_j$, and $g_jyg_j^{-1}$ induces an element of order $p$ on $B_1$ which preserves $\tau^i\circ c_1$ for some $i$, and hence $c_1$, again a contradiction. Finally, the vector of class sizes of $c$ is $\sum_i u_i\,\tau^{i-1}\mathbf n$, which is a cyclic convolution of $\mathbf u$ and $\mathbf n$. For a real vector $\mathbf v$ of length three we have $v_1+v_2\omega+v_3\omega^2=0$, where $\omega$ is a primitive third root of unity, only if $\mathbf v$ is constant. As the discrete Fourier transform turns convolutions into products, the vector of class sizes of $c$ is not constant. So $c$ is good.

\textbf{Step 2.} The case that $\bar X$ is primitive.

If $p\nmid|\bar X|$, a constant colouring is good. If $O_{p'}(\bar X)=1$, then $\bar X$ is a primitive $p$-group, so that $\bar X\cong C_p$ and $n=p$, and we give one point a colour different from the others. So we may assume that $p$ divides $|\bar X|$ and that $O_{p'}(\bar X)\ne 1$. Then $O_{p'}(\bar X)$ is transitive, $p\nmid n$ and $n\ge 4$. If $n=4$, then $p=3$ and $\bar X=A_4$, and the colouring with class sizes $(2,1,1)$ is good, because every element of order $3$ permutes three points cyclically. Let $n\ge 5$. Then $\bar X$ does not contain $\Alt(\Omega)$, since otherwise $\Alt(\Omega)\le O_{p'}(\bar X)$, which forces $p>n$.

Suppose that there is a subset $A$ of $\Omega$ whose setwise stabilizer in $\bar X$ is a $p'$-group. Then $A\ne\emptyset$, and for $a\in A$ the colouring with the classes $A\setminus\{a\}$, $\Omega\setminus A$ and $\{a\}$ is good, since its stabilizer is contained in that of $A$ and $n\ge 5$. By a theorem of Seress \cite{Se}, a primitive group of degree $n>32$ which does not contain $\Alt(\Omega)$ has a regular orbit on the power set of $\Omega$, and so such a set $A$ exists if $n>32$. For $5\le n\le 32$ we used the library of primitive groups in \textsf{GAP} \cite{GAP}; see Appendix A. It turns out that a set $A$ as above exists in all cases except when $n=8$, $p=3$ and $\bar X=\mathrm{A\Gamma L}_1(8)$. In this group there are $56$ elements of order $3$, and $4788$ of the $3^8$ colourings are not preserved by any of them; they satisfy (ii) because $3\nmid 8$.
\end{proof}

\begin{lemma}\label{pgroupset}
Let $Q\le\Sym(\Sigma)$ be a transitive $p$-group, where $p$ is odd. Then there is a subset of $\Sigma$ whose setwise stabilizer in $Q$ is trivial.
\end{lemma}

\begin{proof}
This follows by the argument of Step 1 of the previous proof with two colours: if $U_1\subseteq B_1$ and $U'\subseteq\Sigma$ have trivial stabilizers in the respective groups, put a translate of $U_1$ into the blocks belonging to $U'$ and a translate of $B_1\setminus U_1$ into the other blocks. Since $|B_1|$ is odd, $|U_1|\ne|B_1\setminus U_1|$, and the argument above shows that no element of order $p$ stabilizes the resulting set. The induction starts with $Q\cong C_p$, where a single point has trivial stabilizer.
\end{proof}

\begin{proof}[Proof of Theorem \ref{layer}]
We argue by induction on $m$.

\textbf{Step 1.} We may assume that $X$ is transitive on $\{S_1,\dots,S_m\}$.

Let $E=E_1\times\dots\times E_t$, where the $E_j$ are the products of the $X$-orbits on the set of the $S_i$, and suppose that $t>1$. Put $M_0=X$ and $M_j=C_X(E_1\times\dots\times E_j)$. Then $M_{j-1}/M_j$ is a $p$-nilpotent group which acts faithfully on $E_j$ and contains $\Inn(E_j)$. By induction there are $e_j\in E_j$ such that $C_{M_{j-1}/M_j}(e_j)$ is a $p'$-group, and by Lemma \ref{product} the centralizer of $e_1\cdots e_t$ in $X$ is a $p'$-group.

So let $X$ be transitive on the $S_i$, let $\bar X\le\Sym(m)$ be the induced permutation group and let $B$ be the kernel of this action. We identify all $S_i$ with $S=S_1$.

\textbf{Step 2.} The case that $Y_1$ has at least three orbits on the good elements of $S_1$.

Choose three types, and let $c$ be a colouring of $\{1,\dots,m\}$ as in Lemma \ref{colouring}, applied to $\bar X$. Choose $e_i\in S_i$ good and of the type given by $c(i)$, and let $e=e_1\cdots e_m$. Let $y\in C_X(e)$ be of order $p$. If $y\in B$, then $y$ induces an automorphism of order $p$ on some $S_i$, and this automorphism lies in $Y_i$ and fixes $e_i$, which is not possible. So $y$ induces a permutation $\pi$ of order $p$, and $e_{i\pi}=e_i^{\,y}$ for all $i$. Thus $e_i$ and $e_{i\pi}$ have the same type, so that $\pi$ preserves $c$, a contradiction. Hence $C_X(e)$ is a $p'$-group.

By Lemma \ref{threetypes}, this completes the proof unless $p=3$, $S=\Sz(8)$, $Y_1=\Aut(S)=S\langle\varphi\rangle$ and $\varphi$ permutes the classes of elements of order $7$, and those of order $13$, cyclically. We assume this from now on, and we label these classes by $7_\kappa$ and $13_\kappa$ with $\kappa\in\mathbb Z/3$, in such a way that $\varphi$ maps $\tau_\kappa$ to $\tau_{\kappa+1}$ for $\tau\in\{7,13\}$.

\textbf{Step 3.} A normal form for $X$.

We have $\Aut(E)=\Aut(S)\wr\Sym(m)$, and $\Aut(E)/S^m\cong(\mathbb Z/3)^m\rtimes\Sym(m)$. For $x\in X$ we write $x=(a_1,\dots,a_m)\pi$ with $a_i\in\Aut(S)$ and $\pi\in\Sym(m)$, and we call the vector $f(x)\in(\mathbb Z/3)^m$ with $a_i\in\varphi^{f(x)_i}S$ the field part of $x$. Let $K=O_{3'}(X)$ and let $\bar K$ be its image in $\bar X$. Since $(K\cap B)/S^m$ is a $3'$-subgroup of the $3$-group $B/S^m$, we have $K\cap B=S^m$. Hence $K/S^m$ is a complement of $(\mathbb Z/3)^m$ in $(\mathbb Z/3)^m\rtimes\bar K$, and by the Schur--Zassenhaus theorem it is conjugate to $\bar K$ by an element of $(\mathbb Z/3)^m$. Replacing $X$ by a conjugate under the group $\langle\varphi\rangle^m\le\Aut(E)$, which does not affect the assertion, we may therefore assume that $f(k)=0$ for all $k\in K$. As $X$ normalizes $K$, a short computation in $(\mathbb Z/3)^m\rtimes\Sym(m)$ shows that then $f(x)$ is invariant under $\bar K$ for every $x\in X$. In other words, the field part of every element of $X$ is constant on each orbit $\Delta_j$ of $\bar K$.

\textbf{Step 4.} Construction of $e$.

Let $\chi=(\tau,\kappa)$ be a map from $\{1,\dots,m\}$ to $\{7,13\}\times\mathbb Z/3$, choose $e_i$ in the class $\tau(i)_{\kappa(i)}$ of $S_i$, and let $e=e_1\cdots e_m$. Let $x=(a_1,\dots,a_m)\pi\in C_X(e)$ be of order $3$. If $\pi=1$, then each $a_i$ lies in $C_{\Aut(S)}(e_i)$, which is a $3'$-group, and $a_i^3=1$, so that $x=1$. Hence $\pi$ has order $3$, and $e_{i\pi}$ is the image of $e_i$ under $a_i$. Therefore $\tau(i\pi)=\tau(i)$ and $\kappa(i\pi)=\kappa(i)+f_j$ for $i\in\Delta_j$, where $f_j$ is the value of $f(x)$ on $\Delta_j$. Suppose now that
\begin{equation}\label{majority}
\kappa(i)=0 \text{ for more than half of the elements } i \text{ of each orbit } \Delta_j.
\end{equation}
As $\pi$ permutes the orbits of $\bar K$, comparing the most frequent values of $\kappa$ on $\Delta_j$ and on $\Delta_j\pi$ yields $f_j=0$ for all $j$. Thus $\pi$ preserves $\chi$. So it suffices to find $\chi$ satisfying \eqref{majority} which is not preserved by any element of order $3$ of $\bar X$.

Let $\Sigma$ be the set of the orbits $\Delta_j$. The group induced by $\bar X$ on $\Sigma$ is a transitive $3$-group, and the stabilizer of $\Delta_1$ induces a transitive $3$-nilpotent group $Y$ on $\Delta_1$. If $|\Delta_1|\le 2$, then $Y$ is a $3'$-group, and we let $d$ and $d'$ be the constant maps on $\Delta_1$ with the values $(7,0)$ and $(13,0)$. If $|\Delta_1|\ge 3$, Lemma \ref{colouring} yields a partition of $\Delta_1$ into three sets which is not preserved by any element of order $3$ of $Y$. Refining a partition makes its stabilizer smaller, so, splitting off a point from a class if necessary, we may write the partition as $A\cup B\cup C$ with $1\le|C|<|\Delta_1|/2$. Let $d$ be the map with the values $(7,0)$, $(13,0)$, $(7,1)$ on $A$, $B$, $C$, and let $d'$ be the map with the values $(7,0)$, $(13,0)$, $(13,1)$. In both cases $d$ and $d'$ satisfy \eqref{majority} on $\Delta_1$, they are not preserved by any element of order $3$ of $Y$, and they take their values with different frequencies. By Lemma \ref{pgroupset} there is a subset $U$ of $\Sigma$ with trivial stabilizer in the group induced on $\Sigma$. Fix $g_j\in\bar X$ with $\Delta_1^{g_j}=\Delta_j$, and define $\chi$ on $\Delta_j$ as the translate of $d$ by $g_j$ if $\Delta_j\in U$, and as the translate of $d'$ otherwise. Then $\chi$ satisfies \eqref{majority}, and the argument of Step 1 of the proof of Lemma \ref{colouring} shows that no element of order $3$ of $\bar X$ preserves $\chi$. This completes the proof.
\end{proof}

\section{The quasiprimitive case: structure}

In this section let $p$ be an odd prime and let $(G,V)\in\X$ such that $V$ is quasiprimitive. We put $K=O_{p'}(G)$ and $Z=O_p(G)>1$, and we let $r$ be the characteristic of $V$.

\begin{lemma}\label{structure}
The following hold.
\begin{enumerate}[\rm (a)]
\item Every normal abelian subgroup of $G$ is cyclic and acts fixed point freely on $V-\{0\}$.
\item $Z$ is cyclic.
\item $C_G(K)=Z(K)\times Z$. In particular, every $p$-element of $G$ which is not contained in $Z$ acts nontrivially on $K$ and on $\Fs(K)$, and $C_G(\Fs(K))=Z(\Fs(K))\times Z$.
\item Every component $L$ of $G$ is contained in $K$, and $Z(L)$ is cyclic. A $p$-element of $G$ which normalizes $L$ induces on $L$ the identity or an automorphism which is not inner; in the latter case $L/Z(L)$ is of Lie type over a field $\GF(q_0^p)$.
\end{enumerate}
\end{lemma}

\begin{proof}
(a) If $A\unlhd G$ is abelian, then $V_A$ is homogeneous and faithful, and the assertion follows.

(b) By (a), every characteristic abelian subgroup of $Z$ is cyclic. By a theorem of P.~Hall (see \cite[III, Satz 13.10]{Hu}), $Z$ is the central product of $Z(Z)$, which is cyclic, and a group $E$ which is trivial or extraspecial of exponent $p$. Suppose that $E\ne 1$. Since $K$ centralizes $Z$, the group $G$ acts on $Z/Z(Z)\ne 1$ as a $p$-group, and hence it fixes a nontrivial element $xZ(Z)$, where we may choose $x\in E-Z(E)$. Then $\langle x\rangle Z(Z)$ is an abelian normal subgroup of $G$, hence cyclic, and as $x$ has order $p$, it follows that $x\in Z(Z)\cap E=Z(E)$, a contradiction.

(c) The group $C_G(K)$ is $p$-nilpotent with normal $p$-complement $Z(K)$, and a Sylow $p$-subgroup of $C_G(K)$ centralizes $Z(K)$. Hence it is a characteristic direct factor of $C_G(K)$ and thus contained in $Z$. The remaining statements follow from Lemma \ref{fstar}, applied to the automorphism of $K$ induced by a power of the given element.

(d) We have $L\le K$ as in Step 1 of the proof of Lemma \ref{reduction}. Since $Z(L)$ is an abelian normal subgroup of the normal subgroup $E(K)$ of $G$, and since $V$ is quasiprimitive, a standard argument shows that $Z(E(K))$ is cyclic. The last statement holds since $\Inn(L)$ is a $p'$-group, together with Lemma \ref{out}.
\end{proof}

We now introduce the field which will be used in the remainder of the proof. Let $M=\Fs(K)Z$. Then $M\unlhd G$, so that $V_M$ is homogeneous, and the subalgebra $\mathcal A$ of $\End(V)$ spanned by $M$ over $\GF(r)$ is simple. We let
\[ \F=Z(\mathcal A),\qquad B=C_G(\F),\qquad U=Z(\Fs(K)). \]
Then $\F$ is a field which is normalized by $G$, and $\mathcal A$ is a full matrix algebra over $\F$. Since $U$ and $Z$ are contained in $Z(M)$, they consist of scalars of $\F$; in particular, $p$ divides $|\F|-1$. We may write
\begin{equation}\label{tensor}
V=X\otimes_{\F}Y,
\end{equation}
where $X$ is the simple $\mathcal A$-module, so that $\mathcal A=\End_{\F}(X)\otimes 1$ and $X$ is an absolutely irreducible $\F\Fs(K)$-module. We put $P=\dim_{\F}X$, so that $|V|\ge|\F|^P$. The group $G/B$ is isomorphic to a subgroup of $\Gal(\F/\GF(r))$ and hence cyclic; the elements of $B$ are $\F$-linear, and the elements of $G-B$ act semilinearly on $V$. By Lemma \ref{structure}(c), $C_B(\Fs(K))=U\times Z$, and so
\[ \bar B=B/UZ \]
is isomorphic to a group of automorphisms of $\Fs(K)$.

Let $y\in B$. Since $y$ normalizes $\mathcal A=\End_{\F}(X)\otimes 1$ and is $\F$-linear, the Skolem--Noether theorem shows that
\begin{equation}\label{sn}
y=a_X\otimes a_Y\quad\text{with } a_X\in\GL(X),\ a_Y\in\GL(Y),
\end{equation}
where $a_X$ normalizes the image of $\Fs(K)$ in $\GL(X)$ and induces on it the same automorphism as $y$. If $y$ is a $p$-element which is not a scalar, then $y\notin Z$, and so $a_X$ is not a scalar by Lemma \ref{structure}(c). Write $\Fs(K)=F(K)E(K)$. As $F(K)$ and $E(K)$ centralize each other, $X=X_F\otimes X_E$ for absolutely irreducible modules $X_F$ and $X_E$ of $F(K)$ and $E(K)$, and since $a_X$ normalizes both subgroups, $a_X=a_F\otimes a_E$. Similarly, $X_F$ is the tensor product of modules $X_\ell$ for the Sylow subgroups $O_\ell(K)$ of $F(K)$, and $X_E$ is a tensor product of modules for the components of $K$.

For a semisimple linear map $a$ of a vector space $W$ we denote by $m(a)$ the largest dimension of an eigenspace of $a$ on $W$ (over an algebraic closure). If $a=a_1\otimes a_2$ on $W_1\otimes W_2$, then clearly
\begin{equation}\label{mult}
m(a)\le m(a_1)\dim W_2.
\end{equation}

The tensor factors of $X$ are a priori only defined over an algebraic closure $\bar\F$ of $\F$. The following observation shows that they can be realized over $\F$; it will be used to compare the orders of the groups involved with $|\F|$.

\begin{lemma}\label{realize}
Let $N$ be a normal subgroup of $\Fs(K)$ which is normal in $G$, for instance $N=O_\ell(K)$ or $N$ the product of the components of a given isomorphism type, or let $N$ be a single component of $K$. Let $X_N$ be an irreducible constituent of the $\bar\F N$-module $X\otimes_\F\bar\F$. Then $X_N$ is realizable over $\F$. In particular, if $L$ is a component of $K$ with module $X_L$ of dimension $d$, then $L/Z(L)$ is isomorphic to a subgroup of $\mathrm{PGL}_d(\F)$, and $\F$ contains the field generated by the traces of the elements of $L$ on $X_L$.
\end{lemma}

\begin{proof}
Since $X$ is an absolutely irreducible $\F\Fs(K)$-module and $\Fs(K)$ is a central product of the subgroups considered, $X\otimes\bar\F$ restricted to $N$ is homogeneous, that is, isomorphic to a direct sum of copies of $X_N$. As $X$ is defined over $\F$, the Galois group of $\bar\F/\F$ permutes the irreducible constituents of $X\otimes\bar\F$ restricted to $N$, and so $X_N^\sigma\cong X_N$ for all $\sigma\in\Gal(\bar\F/\F)$. Hence the traces of the elements of $N$ on $X_N$ lie in $\F$, and since Schur indices over finite fields are trivial, $X_N$ is realizable over $\F$. The remaining statements follow, as $L$ acts on $X_L$ faithfully modulo scalars.
\end{proof}

\section{Field automorphisms}

In this section we prove Theorem D. We work in the following situation, which is more general than that of Section 6; in particular, $V$ need not be quasiprimitive.

\begin{hypothesis}\label{galhyp}
Let $p$ be a prime, let $(G,V)\in\X$, and let $\F$ be a subfield of $\End(V)$ which is normalized by $G$ and such that $p$ divides $|\F|-1$. Let $B=C_G(\F)$ and let $\Zh$ be the Sylow $p$-subgroup of $\F^\times$. We put
\[ \Omega=\{0\ne v\in V\mid C_B(v) \text{ is a } p'\text{-group}\}\quad\text{and}\quad \Oh=\{0\ne v\in V\mid C_{B\Zh}(v) \text{ is a } p'\text{-group}\}. \]
\end{hypothesis}

Note that $\Zh$ acts on $V$ as a group of scalars of $\F$, that it centralizes $B$ and is normalized by $G$, but that it need not be contained in $G$. Since $G\cap\Zh$ is a normal $p$-subgroup of $G$, we have $G\cap\Zh\le O_p(G)$. A vector $v\ne 0$ belongs to $\Oh$ if and only if $v$ is not an eigenvector of any $p$-element of $B$ which is not a scalar of $\F$: indeed, if $b\in B$ and $z\in\Zh$ are such that $bz$ is a nontrivial $p$-element fixing $v$, then $b$ is a $p$-element with $bv=z^{-1}v$, and $b$ is not a scalar because scalars have no fixed points; the converse is clear, as the eigenvalues in $\F$ of a $p$-element lie in $\Zh$. The sets $\Omega$ and $\Oh$ are invariant under $G$ and under $\Zh$, and $\Oh\subseteq\Omega$. The group $G/B$ is cyclic, and for $h\in G$ we denote by $\sigma_h$ the automorphism of $\F$ induced by $h$, so that $h(\lambda v)=\sigma_h(\lambda)h(v)$ for $\lambda\in\F$ and $v\in V$.

\begin{theorem}\label{galois}
Assume Hypothesis \ref{galhyp}, and suppose that $\Oh\ne\emptyset$. Then for every normal subgroup $N$ of $G$ we have $\mu_p(N,V)\le|N|_p^{1/2}$, and the action of $N$ on $V$ is not rigid. That is, Theorem C holds for $(G,V)$.
\end{theorem}

\begin{corollary}\label{scalars}
Assume Hypothesis \ref{galhyp}. If every $p$-element of $B$ is a scalar of $\F$, then Theorem C holds for $(G,V)$.
\end{corollary}

For the proof we fix $N\unlhd G$ and write
\[ |N\cap B|_p=p^\gamma\quad\text{and}\quad|NB/B|_p=p^\beta, \]
so that $|N|_p=p^{\beta+\gamma}$. For $0\le j\le\beta$ let $D_j$ be the subgroup of order $p^j$ of the cyclic group $NB/B$.

\begin{lemma}\label{gal1}
Assume Hypothesis \ref{galhyp}.
\begin{enumerate}[\rm (a)]
\item If $v\in\Omega$, then a Sylow $p$-subgroup of $C_N(v)$ is cyclic, intersects $B$ trivially, and has order at most $p^\beta$.
\item If $\sigma$ is an automorphism of $\F$ of order $p^i$, then $[\Zh,\sigma]=\{\sigma(z)z^{-1}\mid z\in\Zh\}$ has order at least $p^i$, and $|\Zh|\ge p^{i+1}$.
\item Let $v\in\Oh$, $b\in N\cap B$ and $u\in\Zh$ such that $bv=uv$. Then $u\in N\cap B\cap\Zh$. In particular, $u$ lies in a subgroup of order at most $p^\gamma$.
\end{enumerate}
\end{lemma}

\begin{proof}
(a) A Sylow $p$-subgroup $P_v$ of $C_N(v)$ satisfies $P_v\cap B\le C_B(v)$, so that $P_v\cap B=1$ and $P_v$ is isomorphic to a subgroup of $NB/B$.

(b) Let $Q_0$ be the order of the fixed field of $\sigma$, so that $|\F|=Q_0^{p^i}$. By Lemma \ref{arith}, $p$ divides $Q_0-1$ and $p^i$ divides $(|\F|-1)/(Q_0-1)$. Since the map $z\mapsto\sigma(z)z^{-1}$ is a homomorphism of the cyclic group $\Zh$ whose kernel is the Sylow $p$-subgroup of the multiplicative group of the fixed field, both assertions follow.

(c) The element $u^{-1}b$ of $B\Zh$ fixes $v$, and so its order is not divisible by $p$. Write $b=b_pb_{p'}$ with commuting elements $b_p$ and $b_{p'}$ of $p$-power order and of order prime to $p$, respectively. Then $u^{-1}b_p$ is a $p$-element which commutes with $b_{p'}$, and thus $u^{-1}b_p=1$. Hence $u=b_p\in N$. Moreover $u\in B$, since $\Zh$ centralizes $\F$, and so $u\in N\cap B\cap\Zh$.
\end{proof}

The next lemma is the key step.

\begin{lemma}\label{ghost}
Assume Hypothesis \ref{galhyp}, let $1\le j\le\beta$ and $v\in\Oh$, and let $S$ be the set of all $z\in\Zh$ such that $|C_N(zv)|_p\ge p^j$. Then
\[ |S|\le|\Zh|\,p^{\gamma-j}. \]
\end{lemma}

\begin{proof}
Fix $g_0\in N$ such that $g_0B$ generates $D_j$, and put $\sigma=\sigma_{g_0}$, which has order $p^j$. Let $z\in S$. As $zv\in\Oh\subseteq\Omega$, Lemma \ref{gal1}(a) shows that $C_N(zv)$ contains a cyclic $p$-subgroup whose image in $NB/B$ contains $D_j$. Hence there is $b_z\in N\cap B$ such that $b_zg_0$ fixes $zv$. Since $b_zg_0(zv)=\sigma(z)\,b_zg_0(v)$, we obtain
\[ b_zg_0(v)=\frac{z}{\sigma(z)}\,v. \]
Now let $z,z'\in S$, and put $b=b_{z'}b_z^{-1}\in N\cap B$ and $w=z'z^{-1}$. Then the last equation, applied to $z$ and to $z'$, gives $bv=uv$ with $u=w\,\sigma(w)^{-1}\in\Zh$. By Lemma \ref{gal1}(c), $u\in T=N\cap B\cap\Zh$, and $|T|\le p^\gamma$. So the homomorphism $\theta:\Zh\to\Zh$, $w\mapsto\sigma(w)w^{-1}$, maps $z'z^{-1}$ into $T$ for all $z,z'\in S$, which means that $S$ is contained in a coset of $\theta^{-1}(T)$. By Lemma \ref{gal1}(b), $|\theta(\Zh)|\ge p^j$, and therefore $|\theta^{-1}(T)|\le|\Zh|\,|T|/p^j$.
\end{proof}

\begin{proof}[Proof of Theorem \ref{galois}]
Fix $v\in\Oh$.

\textbf{Step 1.} $\mu_p(N,V)\le|N|_p^{1/2}$.

If $\gamma\ge\beta$, then $|C_N(v)|_p\le p^\beta\le p^{(\beta+\gamma)/2}$ by Lemma \ref{gal1}(a). So let $\gamma<\beta$ and put $j=\lfloor(\beta+\gamma)/2\rfloor+1$. Then $\gamma<j\le\beta$, and Lemma \ref{ghost} yields $|S|\le|\Zh|p^{\gamma-j}\le|\Zh|/p$. Hence there exists $z\in\Zh-S$, and $|C_N(zv)|_p\le p^{j-1}\le|N|_p^{1/2}$.

For the rest of the proof we assume that the action of $N$ on $V$ is tight and rigid, and we derive a contradiction. Let $|N|_p=p^{2a}$, so that $\mu_p(N,V)=p^a$ and $a\ge 1$, and let $O$ be the $N$-orbit consisting of the minimizers.

\textbf{Step 2.} The case $\gamma\ge\beta$.

For every $w\in\Omega$ we have $p^a\le|C_N(w)|_p\le p^\beta\le p^{(\beta+\gamma)/2}=p^a$. Thus $\beta=\gamma=a\ge 1$ and $\Omega\subseteq O$. We claim that $O$ is invariant under $\Zh$. Indeed, $O=Nv$, and for $n\in N$ and $z\in\Zh$ we have $zn(v)=n(z'v)$ with $z'=n^{-1}zn\in\Zh$, and $z'v\in\Omega\subseteq O$. Since $\Zh$ acts fixed point freely on $V-\{0\}$, it follows that $|\Zh|$ divides $|O|$. But $|\Zh|\ge p^{\beta+1}$ by Lemma \ref{gal1}(b), whereas $|O|_p=|N:C_N(v)|_p=p^{2a-a}=p^\beta$. This is a contradiction.

\textbf{Step 3.} The case $\gamma<\beta$.

Here $j=a+1$ and $a>\gamma$. By Step 1, the $\Zh$-orbit of $v$ contains a minimizer, and so we may assume that $v$ is a minimizer. Let $\langle h\rangle$ be a Sylow $p$-subgroup of $C_N(v)$. By Lemma \ref{gal1}(a), $\langle h\rangle\cap B=1$, and so $\sigma_h$ has order $p^a$. Let
\[ M=\{z\in\Zh\mid zv\in Nv\}=\Zh\cap N\,C_{N\Zh}(v), \]
which is a subgroup of $\Zh$. By Lemma \ref{ghost}, at least $(1-1/p)|\Zh|$ of the vectors $zv$, $z\in\Zh$, satisfy $|C_N(zv)|_p\le p^a$; these vectors are minimizers and thus lie in $O=Nv$. So $|M|>|\Zh|/2$, and hence $M=\Zh$.

Let $z$ be a generator of $\Zh$ and let $n\in N$ with $nv=zv$. As $NB/B$ is abelian, $hnh^{-1}=nc$ for some $c\in N\cap B$. Now $hnh^{-1}(v)=hn(v)=h(zv)=\sigma_h(z)v$. On the other hand, $nv=zv$ implies $n^{-1}v=\sigma_n^{-1}(z^{-1})v$, and thus $nc(v)=\sigma_h(z)v$ gives
\[ cv=\sigma_n^{-1}\bigl(\sigma_h(z)z^{-1}\bigr)v. \]
By Lemma \ref{gal1}(c), the scalar $u=\sigma_n^{-1}(\sigma_h(z)z^{-1})$ lies in $N\cap B\cap\Zh$, which has order at most $p^\gamma$. But $\sigma_h(z)z^{-1}$ generates $[\Zh,\sigma_h]$, and so the order of $u$ is at least $p^a$ by Lemma \ref{gal1}(b). Thus $p^a\le p^\gamma<p^a$, the final contradiction.
\end{proof}

Corollary \ref{scalars} follows from Theorem \ref{galois}, since under its hypothesis $\Oh=V-\{0\}$.

\begin{remark}\label{recover}
If $G$ is a subgroup of the semilinear group $\Gamma(V)$ as in \cite[Proposition 6.3]{K}, we may take $\F=\GF(|V|)$, and then $B$ consists of scalars. So Corollary \ref{scalars} contains \cite[Lemma 6.4]{K}. In the same way one sees that Theorem \ref{galois} contains the arguments of \cite[Lemma 6.14 and Proposition 6.15]{K}.
\end{remark}

\section{Eigenspaces of linear $p$-elements}

We return to the situation and the notation of Section 6. By Theorem \ref{galois}, it remains to find a nonzero vector which is not an eigenvector of any non-scalar $p$-element $y$ of $B$. For this we may assume that $y^p$ is a scalar, and we write $y=a_X\otimes a_Y$ as in \eqref{sn}. In this section we bound the dimensions of the eigenspaces of the tensor factors of $a_X$.

\begin{lemma}\label{extra}
Let $\ell$ be a prime such that $O_\ell(K)$ is not abelian, let $X_\ell$ be the corresponding tensor factor of $X$, of dimension $e=\ell^n$, and let $a\in\GL(X_\ell)$ be an element which normalizes the image of $O_\ell(K)$, is not a scalar, and such that $a^p$ is a scalar. Then there is an integer $s\ge 1$ with $\ell^s\equiv\pm 1\pmod p$ such that
\[ m(a)\le\frac ep\Bigl(1+\frac{p-1}{\ell^s}\Bigr)\le\frac{2e}{p}. \]
Moreover, if $p=3$, $\ell=2$ and $s=1$, then $m(a)=e/2$.
\end{lemma}

\begin{proof}
Note that $\ell\ne p$, and $\ell\ne r$ because $O_r(K)=1$. By Lemma \ref{structure}(a) and the theorem of P.~Hall quoted above, $O_\ell(K)=E\circ C$, where $E$ is extraspecial of order $\ell^{2n+1}$ and $C$ is cyclic, or $\ell=2$ and $C$ is dihedral, semidihedral or generalized quaternion \cite[III, Satz 13.10(b)]{Hu}; a factor $C$ of order $8$ is absorbed into $E$, since then $E\circ C$ is again extraspecial, so that we may assume $|C|\ge16$ in these cases. In the latter cases $\Aut(C)$ is a $2$-group, so that conjugation by $a$ centralizes $C$, and as the corresponding tensor factor of $X_\ell$ is absolutely irreducible, $a$ acts on it as a scalar. Hence we may assume that $C=Z(O_\ell(K))$ is cyclic. Let $\bar E=O_\ell(K)/C$, which is a symplectic space of order $\ell^{2n}$ over $\GF(\ell)$, and let $\ell^{2s}=|[\bar E,a]|$. Since $a$ is not a scalar, $a$ does not centralize $O_\ell(K)$, and so $s\ge 1$. As $\langle a\rangle$ acts on $[\bar E,a]$ without nontrivial fixed points, $\ell^{2s}\equiv 1\pmod p$, so that $\ell^s\equiv\pm 1\pmod p$ and $\ell^s\ge p-1$.

After multiplying $a$ by a scalar in an extension field, we may assume that $a^p=1$. Let $\beta_i$ be the trace of $a^i$ on $X_\ell$, lifted to characteristic zero. Here $(\langle O_\ell(K),a\rangle,O_\ell(K),C,\theta,\varphi)$, with $\theta$ the character of $X_\ell$ restricted to $O_\ell(K)$ and $\varphi$ the faithful linear character of $C$ below it, is a character five in the sense of \cite[Section 3]{Is73}, and $U=C\langle a\rangle$ satisfies $UO_\ell(K)=\langle O_\ell(K),a\rangle$ and $U\cap O_\ell(K)=C$. As the order of $a^i$ modulo $C$ is prime to $|O_\ell(K):C|$, every $a^i$ is good \cite[remark after Definition 3.1]{Is73}, and \cite[Theorem 3.5(a)]{Is73} yields $|\beta_i|^2=|C_{\bar E}(a^i)|=|C_{\bar E}(a)|=\ell^{2n-2s}$ for $p\nmid i$ (the extension factor of $\varphi$ to $U$ is linear and does not change absolute values). As $\ell\ne r$ and $p\ne r$, the Brauer character of $\langle O_\ell(K),a\rangle$ on $X_\ell$ is an ordinary character, so this applies. The multiplicity of an eigenvalue $\lambda$ of $a$ is $\frac1p\sum_{i=0}^{p-1}\lambda^{-i}\beta_i$, and so $m(a)\le(e+(p-1)\ell^{n-s})/p$.

Finally, let $p=3$, $\ell=2$ and $s=1$. Then $\beta_1$ is an element of $\mathbb Z[\omega]$ of absolute value $2^{n-1}$, and as $2$ is a prime in $\mathbb Z[\omega]$, we have $\beta_1=2^{n-1}u$ for a unit $u$. The multiplicities are $(2^n+2^n\mathrm{Re}(\lambda^{-1}u))/3$ for the three third roots of unity $\lambda$. If $u$ is a third root of unity, one of these numbers is $2^{n+1}/3$, which is not an integer. So $-u$ is a third root of unity, and the multiplicities are $0$, $e/2$ and $e/2$.
\end{proof}

\begin{lemma}\label{perm}
Let $L_1,\dots,L_p$ be components of $K$ with modules $X_1,\dots,X_p$ of dimension $d$, and let $a\in\GL(X_1\otimes\dots\otimes X_p)$ be an element which normalizes the image of $L_1\cdots L_p$, permutes the $L_i$ cyclically, and such that $a^p$ is a scalar. Then
\[ m(a)=\frac{d^p+(p-1)d}{p}. \]
\end{lemma}

\begin{proof}
There are isomorphisms $f_i:X_i\to X_{i+1}$ (indices modulo $p$) such that $a(x_1\otimes\dots\otimes x_p)=f_p(x_p)\otimes f_1(x_1)\otimes\dots\otimes f_{p-1}(x_{p-1})$. Identifying $X_{i+1}$ with $X_i$ via $f_i$ for $i<p$, we may assume that $f_1,\dots,f_{p-1}$ are identity maps. Then $a^p=f_p\otimes\dots\otimes f_p$, so that $f_p$ is a scalar and $a$ is a scalar multiple of the cyclic permutation $\pi$ of the tensor factors. The trace of $\pi^i$ is $d$ for $p\nmid i$, and so the eigenvalue $1$ of $\pi$ has multiplicity $(d^p+(p-1)d)/p$, and the other eigenvalues have multiplicity $(d^p-d)/p$.
\end{proof}

Next we consider a component $L$ of $K$ and an element which induces a field automorphism on $L$.

\begin{lemma}\label{frobenius}
Let $W$ be a vector space over a field of characteristic $r$, let $L\le\GL(W)$, and let $a\in\GL(W)$ normalize $L$ such that $a^p$ is a scalar. Suppose that $x\in L$ has prime order $s\ne r$ and that $a$ normalizes $\langle x\rangle$, but does not centralize it. Then
\[ m(a)\le\frac{\dim W}{p}+\Bigl(1-\frac1p\Bigr)\dim C_W(x). \]
\end{lemma}

\begin{proof}
We have $W=C_W(x)\oplus[W,x]$, and both summands are invariant under $a$. On $[W,x]$ the eigenspaces of $x$ belong to nontrivial characters of $\langle x\rangle$, and since $p$ is a prime, $a$ permutes these eigenspaces in orbits of length $p$. Hence every eigenvalue of $a$ on $[W,x]$ has multiplicity $\dim[W,x]/p$.
\end{proof}

\begin{lemma}\label{stable}
Let $L$ be a component of $K$ such that $S=L/Z(L)$ is of Lie type over $\GF(q)$ with $q=\ell^f=q_0^p$, and let $\alpha$ be an automorphism of $L$ of order $p$ which induces a field automorphism on $S$. If $p=3$, assume that $S=\Sz(q)$. Then for every prime $r$ there exist a prime $s\ne r$ and an element $x\in L$ of order $s$ such that $\langle x\rangle$ is $\alpha$-invariant and $C_{\langle x\rangle}(\alpha)=1$.
\end{lemma}

\begin{proof}
Let $k$ and $k'$ be the two largest elements of $I(G)$, and let $r_k$ and $r_{k'}$ be the Zsigmondy primes of Lemma \ref{zsig} (for $p=3$ and $S=\Sz(q)$ we take Zsigmondy primes of $2^{4f}-1$ and $2^f-1$, which divide $q^2+1$ and $q-1$ and do not divide $|\Sz(q_0)|$). These primes are different, and one checks that they do not divide $|Z(L)|$. Let $s\in\{r_k,r_{k'}\}$ with $s\ne r$. By Lemma \ref{complements}(b), $\alpha$ normalizes a Sylow $s$-subgroup $Q$ of $L$, and $C_Q(\alpha)=1$ because $s$ does not divide $|C_S(\alpha)|$ by Lemmas \ref{out} and \ref{zsig}. Thus $\alpha$ acts on the elementary abelian $s$-group $\Omega_1(Z(Q))$. Since $p$ divides $s-1$, this action is diagonalizable over $\GF(s)$, and so there is an $\alpha$-invariant subgroup $\langle x\rangle$ of order $s$.
\end{proof}

\begin{lemma}\label{fieldaut}
Let $L$ be a component of $K$ with module $X_L$ of dimension $d$ over $\F$, where $S=L/Z(L)$ is of Lie type over $\GF(q)$, $q=q_0^p$, and let $a\in\GL(X_L)$ be an element which normalizes the image of $L$, induces a field automorphism $\alpha$ of order $p$ on it, and such that $a^p$ is a scalar. If $p=3$, assume that $S=\Sz(q)$.
\begin{enumerate}[\rm (a)]
\item If $r$ is the defining characteristic of $S$, then $d=d_1^{\,p}$ for an integer $d_1\ge 2$, and $m(a)=(d+(p-1)d_1)/p$.
\item If $r$ is not the defining characteristic of $S$, and if $L$ is generated by $t$ conjugates of the element $x$ of Lemma \ref{stable}, then $m(a)\le d\,(1-(1-1/p)/t)$.
\end{enumerate}
\end{lemma}

\begin{proof}
(a) By Steinberg's twisted tensor product theorem \cite{St} (see \cite[Theorem 5.4.5, Proposition 5.4.6, Remark 5.4.7]{KL}), $X_L\otimes\bar\F$ is a tensor product $M_1\otimes\dots\otimes M_k$ of Galois twists of one module $M$ of dimension $d_M\ge2$, which are permuted by the Frobenius twists. Since $X_L$ is invariant under $\alpha$, which has order $p$ and moves every factor, $\alpha$ permutes the $k$ factors in orbits of length $p$. Grouping the factors accordingly, $X_L\otimes\bar\F\cong N\otimes N^\alpha\otimes\dots\otimes N^{\alpha^{p-1}}$ with $\dim N=d_M^{k/p}$. Now $a$ permutes the $p$ tensor factors $N^{\alpha^i}$ cyclically, and the argument of Lemma \ref{perm} applies.

(b) Since $C_{X_L}(L)=0$, we have $\dim C_{X_L}(x)\le d(1-1/t)$, and the assertion follows from Lemma \ref{frobenius}.
\end{proof}

\section{The linear layer for $p\ge 5$}

In this section let $p\ge 5$, and let $(G,V)\in\X$ with $V$ quasiprimitive. We use the notation of Section 6. Our aim is the following result.

\begin{theorem}\label{linear}
Let $p\ge 5$ and let $(G,V)\in\X$ with $V$ quasiprimitive. Then there exists a nonzero vector $v\in V$ which is not an eigenvector of any $p$-element of $B$ that is not a scalar. In particular, Theorem C holds for $(G,V)$.
\end{theorem}

Note that this holds whether $K$ is solvable or not.

\subsection{Factors and the combination lemma}

We write $X$ as a tensor product of \emph{factors}, which are defined over $\bar\F$ but are realizable over $\F$ by Lemma \ref{realize}: for every prime $\ell$ such that $O_\ell(K)$ is not abelian we have a factor $X_\ell$ of dimension $\ell^{n_\ell}$, and for every isomorphism type of components of $K$ we have one factor, the tensor product of the modules of the $k$ components of this type. Let the factors be $X_c$, $c\in\mathcal C$, with $d_c=\dim X_c\ge 2$, so that $P=\prod_cd_c$ and $|V|\ge|\F|^P$. The group $\bar B=B/UZ$ is isomorphic to a subgroup of the direct product of groups $A_c$, where $A_c$ is the group of automorphisms of $O_\ell(K)$ centralizing $Z(O_\ell(K))$ if $X_c=X_\ell$, and $A_c=\Aut(L)\wr\Sym(k)$ if $X_c$ belongs to $k$ components isomorphic to $L$. Let $\pi_c$ be the projection to $A_c$. We call $c$ \emph{active} if $\pi_c(\bar B)$ contains an element of order $p$, and for active $c$ we let $N_c$ be the number of elements of order $p$ in $\pi_c(\bar B)$ and
\[ \lambda_c=\log_{|\F|}(1+N_c). \]
Moreover, let $\kappa_c$ be a number such that $m(a)\le(1-\kappa_c)d_c$ for every $a\in\GL(X_c)$ which is induced by a $p$-element $y\in B$ with $y^p\in UZ$ and $\pi_c(\bar y)\ne1$. Finally, put $D_c=P/d_c$.

\begin{lemma}\label{combination}
Let $m\ge1$ be the number of active factors. Then the proportion of nonzero vectors of $V$ which are eigenvectors of some non-scalar $p$-element of $B$ is at most
\begin{equation}\label{comb}
\frac{p}{p-1}\sum_{c\ \mathrm{active}}N_c\,|\F|^{-\kappa_cP+\sum_{c'\ne c}\lambda_{c'}},
\end{equation}
where the inner sum runs over the active $c'\ne c$. Suppose moreover that $\lambda_c\le d_c-1$ for all active $c$. If $m=1$, the term of $c$ in \eqref{comb} is $\frac{p}{p-1}N_c|\F|^{-\kappa_cd_cD_c}$, and if $m\ge2$, it is at most
\begin{equation}\label{comb2}
\frac{p}{p-1}N_c\,|\F|^{-D_c(\kappa_cd_c-1)-(m-1)},\qquad\text{where } D_c\ge 2^{m-1}.
\end{equation}
\end{lemma}

\begin{proof}
Let $y\in B$ be a non-scalar $p$-element with an eigenvector $v$ for an eigenvalue in $\F$. Since $U$ and $Z$ consist of scalars, $\bar y\ne1$, and a suitable power $y'$ of $y$ satisfies $|\bar y'|=p$; then $y'$ is not a scalar, $v$ is an eigenvector of $y'$, and $y'^p\in UZ$ is a scalar. The eigenspaces of $y'$ only depend on the subgroup $\langle\bar y'\rangle$, and there are at most $p$ of them. Let $c$ be the first active factor (in a fixed ordering) with $\pi_c(\bar y')\ne1$. By \eqref{mult}, each eigenspace of $y'$ has dimension at most $(1-\kappa_c)\dim_\F V$, and thus at most $|V|\,|\F|^{-\kappa_cP}$ elements. The number of elements of order $p$ of $\bar B$ whose first nontrivial coordinate is $c$ is at most $N_c\prod_{c'\ne c}(1+N_{c'})$, and the number of subgroups of order $p$ is at most $1/(p-1)$ times the number of elements of order $p$. This proves \eqref{comb}.

Now let $\lambda_{c'}\le d_{c'}-1$ for all active $c'$, and let $m\ge2$. Since the dimensions of the factors are at least $2$, the sum of the $d_{c'}$, $c'\ne c$ active, is at most their product, which divides $D_c$. Hence $\sum_{c'\ne c}\lambda_{c'}\le D_c-(m-1)$, and as $\kappa_cP=\kappa_cd_cD_c$, the assertion follows; clearly $D_c\ge2^{m-1}$.
\end{proof}

We say that an active factor $c$ satisfies (PC) if
\begin{equation}\label{pc}
1+N_c\le|\F|^{\kappa_cd_c-1}\qquad\text{and}\qquad\kappa_cd_c\ge 2.
\end{equation}

\begin{corollary}\label{combcor}
If every active factor satisfies \eqref{pc}, then the proportion in Lemma \ref{combination} is less than $1$.
\end{corollary}

\begin{proof}
Note that \eqref{pc} implies $\lambda_c\le d_c-1$, and that $|\F|\ge 8$, since $|\F|\equiv1\pmod p$. If $m=1$, the bound is at most $\frac p{p-1}|\F|^{-1}<1$. If $m\ge2$, then by \eqref{comb2} and \eqref{pc} each term is at most $\frac p{p-1}|\F|^{-(D_c-1)(\kappa_cd_c-1)-(m-1)}\le\frac p{p-1}|\F|^{2-m-2^{m-1}}$, and $\frac p{p-1}\,m\,8^{2-m-2^{m-1}}<1$ for all $m\ge2$.
\end{proof}

If no factor is active, every $p$-element of $B$ is a scalar and Corollary \ref{scalars} applies. It remains to check \eqref{pc} for every type of factor; there will be six exceptional layers, which are treated by explicit computation.

\subsection{Layers}

\begin{lemma}\label{layers}
Let $X_c=X_\ell$ be an active factor of dimension $e=\ell^n$. Then \eqref{pc} holds, unless $(p,\ell,e)$ is one of
\[ (5,2,4),\ (5,2,8),\ (5,2,16),\ (5,2,32),\ (5,3,9),\ (7,2,8). \]
For these six layers $\lambda_c\le e-1$, and the bounds of Lemma \ref{combination} for their terms are as in Table \ref{exclayers}. In particular, the conclusion of Corollary \ref{combcor} also holds if some of the active factors are among these six layers.
\end{lemma}

\begin{table}[h]
\begin{tabular}{c|c|c|c|c}
$(p,\ell,e)$ & $|\F|\ge$ & $\lambda_c$ & term, $m=1$ & term, $m\ge2$\\ \hline
$(5,2,4)$ & $11$ & $1.74$ & $0.061$ & $5.0\cdot10^{-4}$\\
$(5,2,8)$ & $11$ & $2.89$ & $0.0080$ & $5.4\cdot10^{-7}$\\
$(5,2,16)$ & $11$ & $10.19$ & $0.027$ & $4.4\cdot10^{-12}$\\
$(5,2,32)$ & $11$ & $14.93$ & $5.5\cdot10^{-10}$ & $1.7\cdot10^{-30}$\\
$(5,3,9)$ & $16$ & $4.67$ & $0.0020$ & $1.2\cdot10^{-10}$\\
$(7,2,8)$ & $29$ & $4.87$ & $0.026$ & $1.3\cdot10^{-9}$
\end{tabular}
\medskip
\caption{The six exceptional layers. The last column is the bound \eqref{comb2} with $m=2$ and $D_c=2$, which is the largest value for $m\ge 2$.}\label{exclayers}
\end{table}

\begin{proof}
Let $\bar E=O_\ell(K)/Z(O_\ell(K))$ and let $\bar y\in\pi_c(\bar B)$ be of order $p$. Since $p$-elements centralize the dihedral, semidihedral or quaternion part of $O_\ell(K)$ if $\ell=2$ (see the proof of Lemma \ref{extra}), $\bar y$ is determined by its action on $O_\ell(K)/Z$, and the elements of order $p$ of $\pi_c(\bar B)$ lie in $\bar E\rtimes\Sp(2n,\ell)$. Hence $N_c<\ell^{2n}|\Sp(2n,\ell)|<\ell^{2n^2+3n}$. By Lemma \ref{extra}, we may take $\kappa_c=1-2/p$, and $e\ge\ell^s\ge p-1$, where $\ell^{2s}=|[\bar E,\bar y]|$ and $\ell^s\equiv\pm1\pmod p$. In particular $\kappa_ce\ge\frac35\cdot4>2$. Since $U$ and $Z$ contain scalars of order $\ell$ and $p$, we have $|\F|\equiv1\pmod{p\ell}$, and the characteristic of $\F$ is different from $\ell$ and $p$; in particular $|\F|$ is odd if $\ell=2$. Thus \eqref{pc} follows from
\begin{equation}\label{layerineq}
(2n^2+3n)\log_{|\F|}\ell\le\Bigl(1-\frac2p\Bigr)\ell^n-1.
\end{equation}
\emph{The case $\ell\ge5$.} Then $\log_{|\F|}\ell<1$ and $1-2/p\ge3/5$. For $n\ge2$ we have $2n^2+3n\le\frac35\ell^n-1$, and for $n=1$ this holds if $\ell\ge 11$. If $n=1$ and $\ell\in\{5,7\}$, then $\ell\equiv\pm1\pmod p$ is impossible for $p\ge5$, so that the layer is not active.

\emph{The case $\ell=3$.} Then $p\ge5$, $|\F|\ge3p+1\ge16$ and $\log_{|\F|}3\le0.4$, and \eqref{layerineq} holds for $n\ge3$. As $e\ge p-1\ge4$, we have $n\ge2$, and if $n=2$, then $3^s\equiv\pm1\pmod p$ with $s\le2$ forces $p=5$. This gives the exception $(5,3,9)$.

\emph{The case $\ell=2$.} Then $|\F|\ge2p+1\ge11$ and $\log_{|\F|}2\le0.29$, and \eqref{layerineq} holds for $n\ge6$. If $n\le5$, then $2^s\equiv\pm1\pmod p$ for some $s\le5$, so that $p\in\{5,7,11,17,31\}$. For $p\in\{11,17,31\}$ we have $|\F|\ge 2p+1\ge 23$ and $1-2/p\ge9/11$, and \eqref{layerineq} holds for all admissible $n$. For $p=7$ we have $s=3$, $|\F|\ge29$ and $1-2/p=5/7$, and \eqref{layerineq} holds for $n\in\{4,5\}$, but not for $n=3$. For $p=5$ the four layers with $n\in\{2,3,4,5\}$ remain.

We now compute the bounds of Lemma \ref{combination} for the six exceptional layers. For $\bar y$ of order $p$ in $\Sp(2n,\ell)$ with $|[\bar E,\bar y]|=\ell^{2s}$, there are exactly $\ell^{2s}$ elements of order $p$ in the coset $\bar E\bar y$ of $\bar E\langle\bar y\rangle$, as the complements of $\bar E/C_{\bar E}(\bar y)$ are conjugate. The eigenspaces of an element of order $p$ above $\bar y$ have dimension at most $m_s=\lfloor(e+(p-1)\ell^{n-s})/p\rfloor$ by the proof of Lemma \ref{extra}. So the term of $c$ in \eqref{comb} with $m=1$ is at most $\frac p{p-1}\sum\ell^{2s}|\bar y^{\Sp(2n,\ell)}|\,|\F|^{m_s-e}$, the sum running over the classes of elements $\bar y$ of order $p$ of $\Sp(2n,\ell)$, and similarly for \eqref{comb2}. The classes and their centralizers were taken from the character table library of \textsf{GAP} (Appendix A.4): $\Sp(8,2)$ has two classes of elements of order $5$, with centralizers of orders $3600$ ($s=2$) and $300$ ($s=4$), $\Sp(10,2)$ has two such classes, with centralizers of orders $7257600$ ($s=2$) and $1800$ ($s=4$), and $\mathrm{PSp}(4,3)$ and $\Sp(6,2)$ have one class of elements of order $5$ and $7$, respectively, with centralizers of orders $5$ ($s=2$) and $7$ ($s=3$). For $(5,2,4)$ we used that a $5$-nilpotent subgroup of $\Sp(4,2)\cong S_6$ containing an element of order $5$ is cyclic of order $5$, and for $(5,2,8)$ that, by Appendix A.2, the image of $B$ in $\Sp(6,2)$ has a normal $5$-complement $H$ with $|H:C_H(\bar y)|\le16$, so that there are at most $64$ elements of order $5$ in the image. The resulting values are those of Table \ref{exclayers}; the computation is reproduced in Appendix B.

Finally, if some active factors are exceptional layers, then Lemma \ref{combination} still applies, since $\lambda_c\le d_c-1$ for all active factors, and we add the bounds of Table \ref{exclayers} and of the proof of Corollary \ref{combcor} for the other factors. As two exceptional layers have different primes $\ell$, at most two of them occur, and then $p=5$, $\ell\in\{2,3\}$, and $D_c\ge4$ for both. A straightforward computation (Appendix B) shows that the total is at most $0.011$ in all cases with $m\ge2$, and the terms for $m=1$ are less than $0.1$.
\end{proof}

\subsection{Components}

We first collect some consequences of Lemma \ref{realize}. Let $L$ be a component of $K$ with module $X_L$ of dimension $d\ge2$, and $S=L/Z(L)$. Since $X_L$ is realizable over $\F$, we have
\begin{equation}\label{autbound}
|\Aut(S)|\le|S|^2<|\F|^{2d^2},
\end{equation}
where we used that $|\Out(S)|\le|S|$ for every nonabelian simple group $S$.

\begin{lemma}\label{permfactor}
Let $X_c$ be an active factor belonging to $k$ components isomorphic to $L$, and suppose that $L$ does not admit a field automorphism of order $p$ which is induced by an element of $\pi_c(\bar B)$. Then \eqref{pc} holds for $c$.
\end{lemma}

\begin{proof}
Then every element of order $p$ of $\pi_c(\bar B)$ permutes the $k$ components nontrivially, so that $k\ge p\ge5$. By Lemma \ref{perm} and \eqref{mult} we may take $\kappa_c=1-1/p-(p-1)/(pd^{p-1})\ge0.7$, and $d_c=d^k$. Clearly $N_c\le|\Aut(S)|^kk!$.

If $d\ge3$, then by \eqref{autbound}, \eqref{pc} follows from $2kd^2+\log_8(k!)\le0.7d^k-1$, which holds for all $k\ge5$ and $d\ge3$.

Let $d=2$. Then $S\le\mathrm{PGL}_2(\bar\F)$, and by Dickson's theorem either $S\cong\PSL_2(q')$ and $X_L$ is the natural module in the defining characteristic, or $S\cong A_5$ and the characteristic of $\F$ is not $2$ or $5$ (characteristic $3$ is possible, as $\SL_2(5)\le\SL_2(9)$). In the first case $\F$ contains $\GF(q')$, the field of definition of the natural module, so that $|\F|=q'^a$ for some $a$; since $p$ divides $|\F|-1$ but not $q'(q'^2-1)$, the order of $q'$ modulo $p$ is at least $3$ and divides $a$, whence $|\F|\ge q'^3$ and $|\Aut(S)|\le q'^3\log_2q'\le|\F|^{1.6}$. In the second case $p\ge7$, $|\F|$ is odd and $2p$ divides $|\F|-1$, so that $|\F|\ge23$ (the odd prime powers $15<|\F|<23$ are $17$, $19$, which are not $\equiv1$ modulo a prime $p\ge7$) and $|\Aut(S)|=120\le23^{1.6}\le|\F|^{1.6}$. Thus \eqref{pc} follows from $1.6k+\log_8(k!)\le0.7\cdot2^k-1$, which holds for all $k\ge5$. These elementary inequalities are verified in Appendix B.
\end{proof}

Next we consider components with field automorphisms of order $p$. For the cross characteristic case we need generation by few conjugates.

\begin{lemma}\label{tvalues}
Let $S=\PSL_2(q)$ or $S=\Sz(q)$, where $q=q_0^p$, let $L$ be a quasisimple group with $L/Z(L)\cong S$, and let $x\in L$ be an element of prime order $s$, where $s$ is a Zsigmondy prime as in the proof of Lemma \ref{stable}. Then $L$ is generated by two conjugates of $x$.
\end{lemma}

\begin{proof}
It suffices to show that $S$ is generated by two conjugates of the image $\bar x$ of $x$, because $L$ is perfect and a subgroup $J$ with $JZ(L)=L$ satisfies $L=L'\le J$. Note that $\bar x$ has odd prime order $s\ge7$.

If $S=\PSL_2(q)$, then $q\ge32$ and this is \cite[Lemma 3.1]{GS}: an element of odd prime order of $\PSL_2(q)$, $q\ge4$, has $\alpha(\bar x)=2$ unless $q=9$ and $s=3$. (A direct argument via \cite[II, Hauptsatz 8.27]{Hu} must allow for the subfield subgroups $\mathrm{PGL}_2(\sqrt q)$ when $f$ is even, which contain elements of order $s$ dividing $\sqrt q+1$; we therefore cite \cite{GS}.)

Let $S=\Sz(q)$, $q=2^f$ with $f$ odd. Here $s$ is a Zsigmondy prime of $2^f-1$ or of $2^{4f}-1$. By \cite[Theorem 9]{Su}, every proper subgroup of $S$ is isomorphic to $\Sz(q_1)$ with $q=q_1^m$, $m>1$, or is conjugate to a subgroup of the Borel subgroup $H$ of order $q^2(q-1)$, of the dihedral group $B_0$ of order $2(q-1)$, or of one of the groups $B_1,B_2$ of order $4(q\pm\sqrt{2q}+1)$. Since $m$ is odd, $s$ does not divide $|\Sz(q_1)|=q_1^2(q_1^2+1)(q_1-1)$, as the order of $2$ modulo $s$ is $f$ or $4f$ and does not divide $f/m$ or $4f/m$. If $s\mid q-1$, then $s$ divides neither $|B_1|$ nor $|B_2|$, $C_S(\bar x)=T$ is cyclic of order $q-1$, the overgroups of $\bar x$ among the conjugates of $H$ are the two point stabilizers in the doubly transitive action of degree $q^2+1$ of the two fixed points $\omega_1,\omega_2$ of $\bar x$, and the only overgroup among the conjugates of $B_0$ is $N_S(T)$, the stabilizer of $\{\omega_1,\omega_2\}$, whose elements of order $s$ lie in $T$ (compare the proof of \cite[Theorem 12]{Su}). Choose $g\in S$ with $\{\omega_1,\omega_2\}^g\cap\{\omega_1,\omega_2\}=\emptyset$; then $\langle\bar x,\bar x^g\rangle$ lies in no proper subgroup. If $s\mid q\pm\sqrt{2q}+1$, then $s\nmid|H||B_0|$, and $\bar x$ lies in a unique conjugate $N$ of $B_1$ or $B_2$ (the normalizer of the cyclic Hall subgroup containing $\bar x$); as $N$ is not normal, there is $g$ with $\bar x^g\notin N$, and again $\langle\bar x,\bar x^g\rangle=S$.
\end{proof}

\begin{lemma}\label{fieldfactor}
Let $X_c$ be an active factor belonging to $k$ components isomorphic to $L$, where $S=L/Z(L)$ is of Lie type over $\GF(q)$, $q=q_0^p$, and suppose that $\pi_c(\bar B)$ contains an element inducing a field automorphism of order $p$ on some of these components. Then \eqref{pc} holds for $c$.
\end{lemma}

\begin{proof}
\emph{Step 1. Reduction to $k=1$.} Let $\kappa_L$ be the bound for $1-m(a)/d$ given by Lemma \ref{fieldaut} for a single component, and put $\kappa_1=\min(\kappa_L,0.7)$. By Lemmas \ref{perm} and \ref{fieldaut} and \eqref{mult}, we may take $\kappa_c=\kappa_1$ for the block, and $N_c\le|\Aut(S)|^kk!$. Suppose we know that $\log_{|\F|}|\Aut(S)|\le\kappa_1d-1$ (the case $k=1$; note that $N_c\le|\Aut(S)|$ in that case). Then $\kappa_1d\ge2$, and
\[ \kappa_1d^k-1-k(\kappa_1d-1)=\kappa_1(d^k-kd)+k-1\ge2(d^{k-1}-k)+k-1\ge2(2^{k-1}-k)+k-1\ge\log_8(k!), \]
which gives \eqref{pc} for the block, as $|\F|\ge8$. (The last inequality is elementary; see Appendix B.) So it suffices to show $\log_{|\F|}|\Aut(S)|\le\kappa_1d-1$.

\emph{Step 2. The defining characteristic.} By Lemma \ref{fieldaut}(a), $d=d_1^p$ with $d_1\ge2$ and $\kappa_L\ge1-1/p-(p-1)/(pd_1^{p-1})\ge0.75$, so that $\kappa_1=0.7$. Let $\GF(\ell^{f})$ be the smallest field over which $X_L$ is realized, where $q=\ell^e$; by Lemma \ref{realize}, $\GF(\ell^{f})\le\F$. By \cite[Proposition 5.4.6 and Remark 5.4.7]{KL}, $X_L\otimes\bar\F\cong M_1\otimes\dots\otimes M_k$ as in the proof of Lemma \ref{fieldaut}(a), where $k=e/f$ if $S$ is untwisted, $k\ge e/f$ in the twisted cases, and $\dim M_i=d_M\ge2$; in all cases $q\le|\F|^k$. By the proof of Lemma \ref{fieldaut}(a), $p$ divides $k$, so $k\ge5$, and $d=d_M^k$. If $S$ is a Suzuki or Ree group, then \cite[Remark 5.4.7(b)]{KL} gives $d\ge4^k$, and with $|\Aut(S)|\le q^{54}$ it suffices that $54k\le0.7\cdot4^k-1$. Otherwise let $G$ be the corresponding simple algebraic group. If $d_M=2$, then $G$ is of type $A_1$, $S=\PSL_2(q)$ and $|\Aut(S)|\le q^3\log_2q\le q^{3.5}$, as $q\ge32$; so it suffices that $3.5k\le0.7\cdot2^k-1$. If $d_M\ge3$, then $\dim G\le d_M^2$, as $M$ is a nontrivial module for $G$, and $|\Out(S)|\le q^2$; hence $|\Aut(S)|\le q^{d_M^2+2}\le|\F|^{k(d_M^2+2)}$, and it suffices that $k(d_M^2+2)\le0.7d_M^k-1$. All three inequalities hold for $k\ge5$ (Appendix B).

\emph{Step 3. Cross characteristic.} Here $d\ge\delta(S)$, where $\delta(S)$ is the lower bound of Landazuri and Seitz \cite{LS} for the degrees of the nontrivial projective representations of $S$ in cross characteristic, and $q\ge32$, as $q=q_0^p$ with $p\ge5$; moreover $q\ge243$ if $q$ is odd. By Lemma \ref{fieldaut}(b), $\kappa_L\ge(1-1/p)/t\ge0.8/t$, where $t$ is the number of conjugates of the element $x$ of Lemma \ref{stable} needed to generate $L$, and so $\kappa_1\ge0.8/t$ as $t\ge2$. By Lemma \ref{tvalues}, $t=2$ for $\PSL_2(q)$ and $\Sz(q)$. In these two cases we use $|\Aut(S)|\le q^3\log_2q$ and $|\Aut(S)|\le q^5\log_2q$, respectively, together with $|\F|\ge11$ if $p=5$ and $|\F|\ge8$ if $p\ge7$ (in which case $q\ge128$); a direct check (Appendix B) shows that $\log_{|\F|}|\Aut(S)|\le\kappa_1d-1$.

For the remaining groups we use the result of Guralnick and Saxl \cite{GS} that every nontrivial element of $S$ lies in a generating set of at most $t$ of its conjugates, with $t\le n+3$ if $S$ is classical with natural module of dimension $n$, and $t\le 11$ if $S$ is exceptional. Since $|\Out(S)|\le q^2$ for $q\ge32$, we have $|\Aut(S)|\le q^{\dim G+2}$, and as $|\F|\ge 8$ it suffices to show
\begin{equation}\label{ineq}
(\dim G+2)\log_8q\le\frac{0.8\,\delta(S)}{t}-1.
\end{equation}
For the classical groups with $3\le n\le8$ and for all exceptional groups, Table \ref{crosschar} lists both sides of \eqref{ineq} for the smallest admissible even and odd values of $q$, using the bounds $\delta(S)$ of \cite{LS} recalled in Appendix B. The last column of the table shows that \eqref{ineq} holds for every $t$ up to at least $44$, so that the precise values of $t$ from \cite{GS} are not critical. Both sides of \eqref{ineq} are of the form $a\log q$ and $b\,\delta(q)/t-1$, where $\delta(q)$ is a polynomial in $q$ (or in $\sqrt q$) with positive leading coefficient of degree at least $1$, and so \eqref{ineq} remains true for larger $q$ of the same type. Finally, if $S$ is classical with $n\ge9$, then $\delta(S)\ge(q^{n/2-1}-1)/2$ and $\dim G\le n^2$, and \eqref{ineq} holds with $t=2n$ for all $q\ge32$; this is again elementary and checked in Appendix B for $n<200$, the left hand side growing quadratically and the right hand side exponentially in $n$.
\end{proof}

\begin{table}[h]
\small
\begin{tabular}{l|r|r|r|r||l|r|r|r|r}
Type & $q$ & LHS & RHS & max $t$ & Type & $q$ & LHS & RHS & max $t$\\ \hline
$\PSL_3$ & 32 & 16.7 & 135 & 46 & $\PSp_6$ & 243 & 60.8 & $6.4\cdot10^5$ & $9.3\cdot10^4$\\
$\PSU_3$ & 32 & 16.7 & 131 & 45 & $\Omega_7$ & 243 & 60.8 & $2.8\cdot10^8$ & $4.5\cdot10^7$\\
$\PSL_4$ & 32 & 28.3 & 3744 & 894 & $P\Omega_8^\pm$ & 32 & 50.0 & $2.4\cdot10^6$ & $5.3\cdot10^5$\\
$\PSU_4$ & 32 & 28.3 & 3630 & 867 & $G_2$ & 32 & 26.7 & 2380 & 947\\
$\PSp_4$ & 32 & 20.0 & 1756 & 586 & ${}^2G_2$ & 243 & 42.3 & 4276 & 1087\\
$\PSp_4$ & 243 & 31.7 & 3373 & 722 & ${}^3D_4$ & 32 & 50.0 & $2.4\cdot10^6$ & $5.3\cdot10^5$\\
$\PSL_5$, $\PSU_5$ & 32 & 43.3 & $1.0\cdot10^5$ & $1.8\cdot10^4$ & ${}^2F_4$ & 32 & 90.0 & $9.5\cdot10^6$ & $1.1\cdot10^6$\\
$\PSL_6$, $\PSU_6$ & 32 & 61.7 & $2.9\cdot10^6$ & $4.2\cdot10^5$ & $F_4$ & 32 & 90.0 & $1.3\cdot10^{15}$ & $1.5\cdot10^{14}$\\
$\PSL_7$, $\PSU_7$ & 32 & 83.3 & $8.3\cdot10^7$ & $9.9\cdot10^6$ & $E_6$ & 32 & 133 & $2.6\cdot10^{15}$ & $2.1\cdot10^{14}$\\
$\PSL_8$, $\PSU_8$ & 32 & 108 & $2.4\cdot10^9$ & $2.4\cdot10^8$ & ${}^2E_6$ & 32 & 133 & $2.7\cdot10^{21}$ & $2.3\cdot10^{20}$\\
$\PSp_6$, $\PSp_8$ & 32 & $\le63.3$ & $\ge1.4\cdot10^6$ & $\ge3.3\cdot10^5$ & $E_7$ & 32 & 225 & $2.8\cdot10^{24}$ & $1.4\cdot10^{23}$\\
 & & & & & $E_8$ & 32 & 417 & $3.2\cdot10^{42}$ & $8.5\cdot10^{40}$
\end{tabular}
\medskip
\caption{The inequality \eqref{ineq} at the smallest admissible $q$. The columns ``max $t$'' give the largest $t$ for which \eqref{ineq} holds. The rows not listed here (odd $q=243$ for $\PSL_n$, $\PSU_n$, $\PSp_8$, $P\Omega_8^\pm$ and the exceptional groups) all have max $t>1700$; all rows are produced by the script of Appendix B.}\label{crosschar}
\end{table}

\begin{proof}[Proof of Theorem \ref{linear}]
If every $p$-element of $B$ is a scalar, we are done by Corollary \ref{scalars}. Otherwise there are active factors. Active layers are covered by Lemma \ref{layers}, and active factors belonging to components by Lemmas \ref{permfactor} and \ref{fieldfactor}, since by Lemma \ref{structure}(d) an element of order $p$ which normalizes a component induces the identity or a field automorphism on it. Hence, by Corollary \ref{combcor} and Lemma \ref{layers}, there is a nonzero vector which is not an eigenvector of any non-scalar $p$-element of $B$, and Theorem C follows from Theorem \ref{galois}.
\end{proof}

\section{The quasiprimitive case for $p=3$}

In this section let $p=3$, let $(G,V)\in\X$ with $V$ quasiprimitive, and suppose that $K$ is nonsolvable; if $K$ is solvable, then $G$ is solvable and Theorem C holds by \cite[Section 6]{K}. By the classification of finite simple groups, the nonabelian composition factors of $K$ are Suzuki groups $\Sz(q)$ with $q=2^f$, $f\ge3$ odd. We keep the notation of Section 9. We again use Theorem \ref{galois}, but for $p=3$ some factors are too small for \eqref{pc}, and we therefore work with the bound \eqref{comb} directly. Put
\begin{equation}\label{crit3}
\Theta=\frac32\sum_{c\ \mathrm{active}}N_c\,\prod_{c'\ne c}(1+N_{c'})\,|\F|^{-\kappa_cP}.
\end{equation}
By Lemma \ref{combination} and Theorem \ref{galois}, Theorem C holds for $(G,V)$ if $\Theta<1$.

\begin{lemma}\label{data3}
The factors of $X$ and the numbers $d_c$, $N_c$, $\kappa_c$ and $|\F|$ satisfy the following.
\begin{enumerate}[\rm (a)]
\item $|\F|\equiv1\pmod3$, and $|\F|\equiv1\pmod\ell$ for every layer $X_\ell$; the characteristic $r$ of $\F$ is different from $3$ and from all such $\ell$.
\item For a layer $X_\ell$ of dimension $e=\ell^n$, we have $N_c\le\ell^{2n}|\Sp(2n,\ell)|$, and we may take $\kappa_c=1/2$ if $\ell=2$ and $\kappa_c=1-(1+2/\ell)/3$ if $\ell\ge5$.
\item Every component $L$ of $K$ satisfies $L/Z(L)\cong\Sz(q)$, and $|Z(L)|\le2$, with $|Z(L)|=2$ only if $q=8$ and $r\ne2$. If $X_c$ belongs to $k$ components $\Sz(q)$, then $N_c\le|\Aut(\Sz(q))|^kk!$, and if $k=1$, then $N_c\le2|\Sz(q):\Sz(q^{1/3})|$ if $3\mid f$, and $c$ is not active if $3\nmid f$. If $L\cong2.\Sz(8)$, then no element of $G$ induces a nontrivial automorphism of order $3$ on $L$; the module of $L$ has dimension at least $8$, and at least $16$ unless $r=5$; and if $X_c$ belongs to $k$ such components, then $c$ is active only if $k\ge3$, and $N_c\le|\Sz(8)|^kk!$.
\item If $r=2$, the modules of the components are Steinberg tensor products of Frobenius twists of the natural module, of dimensions $4^j$ with $1\le j\le f$, and $\F$ contains the field of definition $\GF(2^{g_0})$ of such a module, where $g_0j\ge f$. In particular $q\le|\F|^j$, and $|\Aut(\Sz(q))|<q^6\le|\F|^{6j}$. A component on which an element of $B$ induces a field automorphism of order $3$ satisfies $3\mid j$, and then $\kappa_c\ge1-(d_1^3+2d_1)/(3d_1^3)$ with $d_1=4^{j/3}$. If $r$ is odd, then $|\F|\ge7$, the module of each component with $Z(L)=1$ has dimension at least $(q-1)\sqrt{q/2}$, and we may take $\kappa_c=1/3$ for field automorphisms.
\item If an element of order $3$ of $\pi_c(\bar B)$ permutes three components with modules of dimension $d$, then we may take $\kappa_c=1-(d^3+2d)/(3d^3)$.
\item If $E(K)=1$, then there is a layer $X_\ell$ with $\ell=2$ and $n\ge6$, or with $\ell$ odd and $\ell^n\ge121$. If $\ell=2$ and $n=6$, then $N_c\le2^{12}\cdot3\cdot7\cdot|\Sz(8)|$.
\end{enumerate}
\end{lemma}

\begin{proof}
(a) is clear since $U$ and $Z$ consist of scalars of $\F$, and $O_r(K)=1$. (b) follows as in the proof of Lemma \ref{layers}, using Lemma \ref{extra} with $p=3$. (c) The Schur multiplier of $\Sz(q)$ is trivial for $q>8$ and elementary abelian of order $4$ for $q=8$, and $Z(L)$ is cyclic and acts by scalars; in characteristic $2$ there are no nontrivial scalars of $2$-power order. The bound for general $k$ is clear, as $\Aut(L)$ embeds into $\Aut(\Sz(q))$. Let $L\cong2.\Sz(8)$. The outer automorphism of order $3$ of $\Sz(8)$ acts fixed point freely on the Schur multiplier (the group $2^2.\Sz(8).3$ has trivial centre, see \cite{ATLAS}), so it permutes the three double covers of $\Sz(8)$ and does not lift to an automorphism of $L$. Hence $\Aut(L)$ is isomorphic to $\Sz(8)$, of order prime to $3$, and an element of order $3$ normalizing $L$ centralizes it. The bounds for the dimension follow from the character table of $2.\Sz(8)$ and its Brauer tables for the primes $5$, $7$ and $13$ dividing $|\Sz(8)|$ \cite{ATLAS,GAP}: the smallest degrees of faithful irreducible representations are $40$ in characteristic $0$ (and hence in every characteristic not dividing $|\Sz(8)|$), $8$ in characteristic $5$, $40$ in characteristic $7$ and $16$ in characteristic $13$ (see Appendix A.3). This is consistent with \cite{LS}, where only the bound $8$ is given in this case. If $k=1$, then $\pi_c(\bar B)\le\Aut(\Sz(q))$, and since $\Sz(q)$ is a $3'$-group whose outer automorphism group is cyclic of order $f$, Lemma \ref{complements}(a) shows that the elements of order $3$ of $\Aut(\Sz(q))$ are the conjugates of $\varphi^{\pm1}$, where $\varphi$ is a field automorphism of order $3$; there are $2|\Sz(q):C_{\Sz(q)}(\varphi)|=2|\Sz(q):\Sz(q^{1/3})|$ of them.

(d) The first statement is Steinberg's tensor product theorem \cite{St} for $\Sz(q)$: the irreducible modules are $\bigotimes_{i\in I}N^{(2^i)}$ for subsets $I$ of $\mathbb Z/f$, of dimension $4^{|I|}$. Such a module is realized over $\GF(2^{g_0})$, where $g_0$ is the smallest positive integer with $I+g_0=I$, and not over a smaller field; since $I$ is a union of cosets of the subgroup of order $f/g_0$, we have $j=|I|\ge f/g_0$. By Lemma \ref{realize}, $\GF(2^{g_0})\le\F$. If an element of $B$ induces a field automorphism $\alpha$ of order $3$ on a component, then the module is $\alpha$-invariant, so that $I+f/3=I$ and $3\mid j$, and the bound for $\kappa_c$ is Lemma \ref{fieldaut}(a). If $r$ is odd, then $|\F|\ge7$ as $|\F|\equiv 1\pmod 3$, the dimension bound is \cite{LS}, and $\kappa_c=1/3$ follows from Lemmas \ref{fieldaut}(b) and \ref{tvalues}.

(e) is Lemma \ref{perm}.

(f) If $E(K)=1$, then $C_K(F(K))\le F(K)$, and $K/F(K)$ is involved in the direct product of the groups $\Sp(2n_\ell,\ell)$ and an abelian group. As $K/F(K)$ has a composition factor $\Sz(q)$, $q=2^f$, some $\Sp(2n_\ell,\ell)$ involves $\Sz(q)$. Let $r_1$ and $r_4$ be Zsigmondy primes of $2^f-1$ and $2^{4f}-1$, respectively; they exist as $f\ge3$ and $4f\ge12$. Both divide $|\Sz(q)|=q^2(q-1)(q^2+1)$, and hence $|\Sp(2n,\ell)|=\ell^{n^2}\prod_{i=1}^n(\ell^{2i}-1)$. Moreover $r_1\equiv1\pmod f$ and $r_4\equiv1\pmod{4f}$, so that $r_1\ge7$ and $r_4\ge13$, and $r_1\ne r_4$.

If $\ell=2$, then $r_4$ divides $2^{2i}-1$ for some $i\le n$, and since $2$ has order $4f$ modulo $r_4$, we get $4f\mid 2i$ and thus $2n\ge4f\ge12$. In particular, if $n=6$, then $f=3$ and $q=8$.

If $\ell$ is odd, then $n\ge2$ by Dickson's theorem, as $\Sp(2,\ell)=\SL_2(\ell)$ does not involve a Suzuki group. If $n=2$ and $\ell\in\{5,7\}$, then $|\Sp(4,5)|=2^7\cdot3^2\cdot5^4\cdot13$ and $|\Sp(4,7)|=2^9\cdot3^2\cdot5^2\cdot7^4$ have at most one prime divisor which is at least $7$, so that they cannot be divisible by both $r_1$ and $r_4$. Hence $\ell^n\ge\min(5^3,7^3,11^2)=121$.

Finally, let $\ell=2$ and $n=6$, so that $q=8$ and $r_4=13$. Since $2$ has order $12$ modulo $13$, an element of order $13$ in the image $\bar G$ of $G$ in $\Sp(12,2)$ acts irreducibly on $\GF(2)^{12}$, and so $\bar G$ is an irreducible subgroup of $\GL_{12}(2)$ containing a $\mathrm{ppd}(12,2;12)$-element. By the Main Theorem of \cite{GPPS}, $\bar G$ is one of the groups of Examples 2.1--2.9 there. Since $\bar G$ is $3$-nilpotent, its nonabelian composition factors are Suzuki groups, and one of them is $\Sz(8)$. This excludes the classical examples 2.1, whose normal subgroups $\SL_d$, $\Sp_d$, $\mathrm{SU}_d$, $\Omega_d$ have order divisible by $3$, as well as the alternating and sporadic examples 2.6 and 2.7. The reducible examples 2.2 are excluded by irreducibility. The imprimitive examples 2.3 require $e<d$, and the symplectic type examples 2.5 require $q$ odd. In the extension field examples 2.4, $\bar G\le\Gamma\mathrm L_{12/b}(2^b)$ with $b\mid12$; for $b=12$ the group is solvable, and for $b\in\{2,4,6\}$ the group $\Sz(8)$ has no nontrivial module of dimension $12/b$ over $\GF(2^b)$, as its nontrivial modules in characteristic $2$ have dimension at least $4$ and need the field $\GF(8)$. So $b=3$, and $\bar G\cap\GL_4(8)$ is one of the examples of \cite{GPPS} in dimension $4$ over $\GF(8)$; as before the only possibility with a composition factor $\Sz(8)$ is ${}^2B_2(8)$ in its natural representation \cite[Table 6]{GPPS}. Finally, $\Sz(8)$ has no absolutely irreducible module of dimension $12$ in characteristic $2$, which excludes the remaining examples 2.8 over $\GF(2)$, and the examples 2.9 concern groups of Lie type in odd characteristic. Hence $\bar G$ normalizes a subgroup $\Sz(8)$ acting on $\GF(8)^4$. Its normalizer in $\Sp(12,2)$ has order dividing $3\cdot7\cdot|\Sz(8)|$, and the bound for $N_c$ follows.
\end{proof}

\begin{lemma}\label{tail3}
If $P\ge2^{10}$, then $\Theta<1$.
\end{lemma}

\begin{proof}
Let $\Lambda=\sum_c\lambda_c$, the sum over all active factors, and $L_P=\log_2P$. All $\kappa_c$ are at least $1/3$, $|\F|\ge4$, and the number $m$ of active factors is at most $L_P$. Hence $\Theta\le\frac32m|\F|^{\Lambda-P/3}$, and it suffices to show that $\Lambda+\log_4(1.5L_P)<P/3$. We claim that $\lambda_c\le H(L_c)$ for every factor, where $L_c=\log_2d_c$ and $H(x)=2x^2+4x+1+\frac x2\log_2x$. For a layer, $1+N_c\le\ell^{2n^2+3n}$ and $|\F|>\ell$ by Lemma \ref{data3}(a), so $\lambda_c\le2n^2+3n\le2L_c^2+3L_c$. For components in characteristic $2$, Lemma \ref{data3}(c),(d) gives $\lambda_c\le6\sum j+\log_4(k!)+1\le3L_c+\frac{L_c}{2}\log_2L_c+1$. In cross characteristic, if $L/Z(L)=\Sz(q)$ and $Z(L)=1$, then $q\le(1.62d)^{2/3}$ and $|\Aut(\Sz(q))|\le q^{6}$, so that $\log_7|\Aut(\Sz(q))|\le 2\log_2d$; if $L\cong2.\Sz(8)$, then $d\ge8$ and $\log_7|\Aut(L)|=\log_7|\Sz(8)|<6=2\log_2 8$. In both cases $\lambda_c\le2L_c+\log_7(k!)+1$. This proves the claim. Since $\sum_cL_c=L_P$ and $m\le L_P$, we get $\Lambda\le2L_P^2+5L_P+L_P\log_2L_P$. For $L_P\ge10$ we have $2^{L_P}/3>2L_P^2+5L_P+L_P\log_2L_P+\log_4(1.5L_P)$, which completes the proof.
\end{proof}

\begin{proposition}\label{p3}
Let $p=3$, let $(G,V)\in\X$ with $V$ quasiprimitive, and suppose that $K$ is nonsolvable. Then $\Theta<1$. In particular, Theorem C holds for $(G,V)$.
\end{proposition}

\begin{proof}
By Lemma \ref{tail3} we may assume that $P<2^{10}$. By Lemma \ref{data3}(c),(f), one of the factors carries a Suzuki group: it is a block of Suzuki components, or a layer as in Lemma \ref{data3}(f). In particular all factors other than this one have dimension less than $2^8$. Thus the possible configurations of factors are: at most one layer $X_\ell$ for each prime $\ell\ne3$, with $\ell^n<2^{10}$; blocks of Suzuki components, of dimension $\prod_i4^{j_i}$ in characteristic $2$, or at most two components of type $\Sz(8)$, $\Sz(32)$ or $\Sz(128)$, or at most three components of type $2.\Sz(8)$, in odd characteristic; and the product of all dimensions is less than $2^{10}$. For each configuration we computed an upper bound for $\Theta$, using Lemma \ref{data3}, the smallest value of $|\F|$ compatible with Lemma \ref{data3}(a),(d), and the smallest dimensions allowed by Lemma \ref{data3}(c),(d) in odd characteristic; for $2.\Sz(8)$ we treated the cases $r=5$ (dimension $8$) and $r\ne5$ (dimension $16$) separately. For components in characteristic $2$ we used the bounds of Lemma \ref{data3}(d) in terms of $|\F|$, which are uniform in $q$, and we assumed that an element inducing a field automorphism exists whenever $3\mid j$; this can only increase the bound. These choices give upper bounds for $\Theta$ for all configurations of the given shape. Indeed, larger dimensions in odd characteristic increase $P$ and leave the $N_c$ unchanged, and the program verifies that in every term of $\Theta$ the total exponent of $|\F|$ is negative, so that the bound decreases as $|\F|$ increases.

The computation (Appendix B.4) covers $411$ distinct configurations, $386$ of which contain active factors, and in all of them $\Theta<1$. An independent implementation (Appendix B.5), with a different enumeration, an independent search for the smallest admissible $|\F|$ and high-precision arithmetic, produces the same set of configurations and the same values. It also tests the monotonicity used above, by recomputing $\Theta$ for the three successive admissible values of $|\F|$ following the smallest one and for larger dimensions in odd characteristic; no increase occurs. The largest values of $\log_{|\F|}\Theta$ are listed in Table \ref{tight3}. The tightest case is a single component $\Sz(8)$ in dimension $14$ over $\GF(7)$, with $\Theta\le7^{-0.36}<0.5$.
\end{proof}

\begin{table}[h]
\begin{tabular}{l|c|c|c}
Factors & $|\F|$ & $P$ & $\log_{|\F|}\Theta\le$\\ \hline
$\Sz(8)$, cross characteristic & $7$ & $14$ & $-0.36$\\
$\Sz(8)$, cross char., and a layer $Q_8$ & $7$ & $28$ & $-3.37$\\
$2.\Sz(8)$ in characteristic $5$, and a layer $Q_8$ & $25$ & $16$ & $-6.89$\\
$\Sz(8)$, cross char., and a layer $2^{1+4}$ & $7$ & $56$ & $-9.55$\\
$\Sz(q)$, natural module, and a layer $5^{1+2}$ & $4096$ & $20$ & $-9.66$\\
$2.\Sz(8)$ in characteristic $5$, and a layer $2^{1+4}$ & $25$ & $32$ & $-12.97$\\
$\Sz(q)$, natural module, and a layer $7^{1+2}$ & $64$ & $28$ & $-13.57$\\
$2.\Sz(8)$ in characteristic $\ne5$, and a layer $Q_8$ & $7$ & $32$ & $-14.16$
\end{tabular}
\medskip
\caption{The configurations with the largest values of $\Theta$ for $p=3$.}\label{tight3}
\end{table}

\section{Proofs of the main results}

\begin{proof}[Proof of Theorem C]
Let $p$ be a prime. If $p=2$, then $G$ is solvable because $O_{p'}(G)$ has odd order, and Theorem C holds by \cite[Theorem C]{K}. So let $p$ be odd. We argue by induction on $|G|$, so that Hypothesis \ref{ind} holds. If $V$ is not quasiprimitive, we are done by Proposition \ref{imprimitive}. If $V$ is quasiprimitive and $p\ge 5$, we are done by Theorem \ref{linear}. If $p=3$ and $V$ is quasiprimitive, we are done by \cite[Section 6]{K} if $K$ is solvable, and by Proposition \ref{p3} otherwise.
\end{proof}

\begin{proof}[Proof of Theorem A]
If $p=2$, then $H$ is solvable, and the assertion is \cite[Theorem A]{K}. So let $p$ be odd, and let $G=HP$ be the semidirect product. Then $G$ is $p$-nilpotent with $O_{p'}(G)=H$, and $P$ acts faithfully on $H$. By Theorem C, the hypothesis of Lemma \ref{reduction} is satisfied. Hence there exists $x\in H$ with $|C_G(x)|_p\le|G|_p^{1/2}=|P|^{1/2}$, and by Lemma \ref{classes} the class $c=x^H$ satisfies $|C_P(c)|=|C_G(x)|_p\le|P|^{1/2}$.
\end{proof}

\begin{proof}[Proof of Corollary B]
Since $P$ is solvable and the action is coprime, every subgroup $Q\le P$ fixes the same number of irreducible characters and of conjugacy classes of $H$, and hence the actions of $P$ on $\Irr(H)$ and on $\cl(H)$ are permutation isomorphic \cite[Theorem 13.24 and Lemma 13.23]{Is} (compare \cite[Section 7]{K}). Thus Corollary B follows from Theorem A.
\end{proof}

\section{A nonsolvable family of examples}

In this section we prove Theorem E. Let $p\equiv1\pmod 3$ be a prime, and let $q_0$ be a prime power such that $q_0$ has order $3$ modulo $p$. Put $q=q_0^p$, and let $\F=\GF(q^3)$, $\F_0=\GF(q_0^3)$ and $V=\F^2$. Let $S=\SL_2(q)$ act on $V$ by matrices over $\GF(q)$, let $W=\F^\times$ act on $V$ by scalar multiplication, and let $Z$ and $U$ be the Sylow $p$-subgroup and the Hall $p'$-subgroup of $W$, respectively. Let $\varphi$ be the map $(x_1,x_2)\mapsto(x_1^{q_0^3},x_2^{q_0^3})$ of $V$, which has order $p$ and fixed field $\F_0$ on scalars. We put $K=SU$, $H_0=SW=K\times Z$ and $G=H_0\langle\varphi\rangle$. Let $p^c=(q_0^3-1)_p$ and $Z_1=\Omega_1(Z)$.

We need the following observation.

\begin{lemma}\label{regular}
Let $d\ge 1$ and let $q$ be a prime power. Consider $V_d=\GF(q^{d+1})^d$ as the space of $d\times(d+1)$-matrices over $\GF(q)$, on which $\GL_d(q)$ acts from the left and $\GF(q^{d+1})^\times$ acts by scalar multiplication. Then $\GL_d(q)\circ\GF(q^{d+1})^\times$ acts regularly on the set of matrices of rank $d$.
\end{lemma}

\begin{proof}
Let $(g,s)$ fix a matrix $v$ of rank $d$, and let $R$ be the $\GF(q)$-span of the coordinates of $v$, which is a hyperplane of $\GF(q^{d+1})$. Then $Rs=R$, so that $R$ is a vector space over $\GF(q)[s]$. The degree of $\GF(q)[s]$ over $\GF(q)$ divides $d$ and $d+1$, and thus $s\in\GF(q)$. Now $gv=s^{-1}v$, and as the rows of $v$ are linearly independent over $\GF(q)$, we get $g=s^{-1}$, so that $(g,s)$ acts trivially. The number of matrices of rank $d$ is $\prod_{i=0}^{d-1}(q^{d+1}-q^i)$, which equals $|\GL_d(q)|(q^{d+1}-1)/(q-1)=|\GL_d(q)\circ\GF(q^{d+1})^\times|$.
\end{proof}

\begin{theorem}\label{examples}
In the situation above, the following hold.
\begin{enumerate}[\rm (a)]
\item $(G,V)\in\X$, $O_{p'}(G)=K$ is nonsolvable, and $|G|_p=p^{c+2}$.
\item $|C_G(v)|_p=p$ for every $0\ne v\in V$.
\item $N=KZ_1\langle\varphi\rangle$ is a normal subgroup of $G$ with $|N|_p=p^2$, and $|C_N(v)|_p=p$ for every $0\ne v\in V$. In particular, the action of $N$ on $V$ is tight, but not rigid.
\end{enumerate}
\end{theorem}

\begin{proof}
(a) Since $q\equiv q_0\pmod p$, the order of $q$ modulo $p$ is $3$, so that $p$ divides $q^3-1$ but not $q(q^2-1)=|S|$. Hence $Z>1$ and $K$ is a nonsolvable $p'$-group. Clearly $K$ is normal in $G$ and $G/K\cong Z\rtimes\langle\varphi\rangle$, so that $G$ is $p$-nilpotent and $Z\le O_p(G)$. The span of $H_0$ in $\End(V)$ is the full matrix algebra $M_2(\F)$, and so $V$ is a faithful irreducible $G$-module. By Lemma \ref{arith}, $(q^3-1)_p=p(q_0^3-1)_p$, so that $|Z|=p^{c+1}$ and $|G|_p=p^{c+2}$.

We next note two facts. First, $C_Z(\varphi)=Z\cap\F_0$ has index $p$ in $Z$, so that $[Z,\varphi]=Z_1$ and $\varphi$ acts trivially on $Z/Z_1$. Second, every $H_0$-orbit on $V-\{0\}$ contains a vector which is fixed by $\varphi$. Indeed, $H_0$ is transitive on the vectors of rank $1$, and $u\otimes1$ with $u\in\GF(q_0)^2$ is fixed by $\varphi$. By Lemma \ref{regular}, $H_0$ acts semiregularly on the vectors of rank $2$, with $\gcd(2,q-1)$ orbits, which are distinguished by the determinant modulo squares. If the coordinates of $v_0$ form a basis over $\GF(q_0)$ of a hyperplane of $\F_0$, then they are linearly independent over $\GF(q)$, as $\F_0$ and $\GF(q)$ are linearly disjoint over $\GF(q_0)$. Thus the vectors $hv_0$ with $h\in\GL_2(q_0)$ have rank $2$, they are fixed by $\varphi$, and they meet every $H_0$-orbit, because $|\GF(q)^\times:\GF(q_0)^\times|$ is odd.

(b) Since $H_0=K\times Z$ and $Z$ acts fixed point freely on $V-\{0\}$, the stabilizers $C_{H_0}(v)$ are $p'$-groups, so that $|C_G(v)|_p\le|G:H_0|=p$. By the second fact, every $H_0$-orbit is $\varphi$-invariant and hence a $G$-orbit. Thus $|C_G(v)|=p|C_{H_0}(v)|$ and $|C_G(v)|_p=p$.

(c) Let $H_1=K\times Z_1$. By the first fact, $G/H_1\cong Z/Z_1\times\langle\varphi\rangle$ is abelian, so that $N=H_1\langle\varphi\rangle$ is normal in $G$ and $|N|_p=p^2$. As in (b), $|C_N(v)|_p\le p$. Let $O$ be an $H_0$-orbit and $v_O\in O$ a vector fixed by $\varphi$. The stabilizers of $H_0$ are contained in $K\le H_1$, so the $H_1$-orbits in $O$ are the sets $H_1zv_O$ with $z$ running over representatives of the cosets of $Z_1$ in $Z$. Now $\varphi$ maps $H_1zv_O$ to $H_1\varphi(z)v_O$, and $\varphi(z)\in zZ_1$ by the first fact. So every $H_1$-orbit is $\varphi$-invariant and hence an $N$-orbit, and as before $|C_N(v)|_p=p=|N|_p^{1/2}$. Thus the action of $N$ is tight, and every nonzero vector is a minimizer. Since the vectors of rank $1$ and those of rank $2$ lie in different $N$-orbits, the action is not rigid.
\end{proof}

Theorem E follows from Theorem \ref{examples}. The smallest example is obtained for $p=7$ and $q_0=2$, so that $V=\GF(2^{21})^2$ and $|G|_7=7^3$.

\begin{remark}\label{examplesremark}
(a) In these examples every $p$-element of $B=H_0$ is a scalar, and so Corollary \ref{scalars} applies. On the other hand, they show that for $p\ge 7$ Theorem C cannot be proved by finding a vector with a $p'$-centralizer in $G$.

(b) More generally, let $d\ge 1$, let $q_0$ have order $d+1$ modulo $p$, and let $q=q_0^p$. Then the group $(\SL_d(q)\circ\GF(q^{d+1})^\times)\rtimes\langle\varphi\rangle$, acting on $\GF(q^{d+1})^d$, has the property that $|C(v)|_p=p$ for all vectors $v$ of rank $d$, by Lemma \ref{regular} and the argument of the proof of Theorem \ref{examples}. For $d=1$ these are groups as in \cite[Theorem B]{K}, and for $d=2$ they are the groups above. Such examples exist whenever $d+1$ divides $p-1$; in particular, they do not exist for $p=3$ and $d\ge 2$.
\end{remark}

\section{Outlook}

We conclude with some remarks.

First, the proof of Theorem A for $p\ge 5$ is uniform in the sense that it does not distinguish between solvable and nonsolvable groups in the quasiprimitive case, and so it gives a proof of \cite[Theorem C]{K} for $p\ge 5$ which does not depend on the classification of Holt and Yang \cite{HY}. For $p=3$ the solvable quasiprimitive case still relies on \cite{HY} and on computer calculations. It would be desirable to extend the counting arguments of Section 9 to $p=3$. The difficulty is that for $p=3$ the eigenspaces of a linear element of order $3$ may have half the dimension of the module, so that Lemma \ref{combination} does not apply to small factors.

Second, Theorem \ref{galois} shows that in quasiprimitive linear actions of $p$-nilpotent groups, the elements acting by field automorphisms never prevent the bound $|N|_p^{1/2}$. It might be interesting to see whether this principle can be used for other questions on orbit sizes of semilinear groups.

Third, the question raised in \cite[Section 8]{K} remains open: if a group $N$ acts faithfully and coprimely on a group $H$, does there always exist a class $c$ of $H$ with $|C_N(c)|\le|N|^{1/2}$? By Theorem A, the answer is affirmative if $N$ is a $p$-group, and by a theorem of Halasi and Podoski \cite{HP}, the analogous statement for elements is true.

Finally, the examples of Theorem E and Remark \ref{examplesremark} suggest the problem of classifying all coprime actions of $p$-groups $P$ on arbitrary finite groups $H$ for which $\min_{c\in\cl(H)}|C_P(c)|=|P|^{1/2}$.

\appendix
\section{Computations with \textsf{GAP}}

All computations were carried out with \textsf{GAP} 4.12.1 \cite{GAP}, using the libraries of primitive groups (which is complete for degrees below $4096$, see \cite{CQRD}), of character tables, and of tables of marks. The input files and the complete output are available as supplementary material; the input files are reproduced below.

\textbf{A.1. The primitive groups of Lemma \ref{colouring}.} For every primitive group $\bar X$ of degree $5\le n\le32$ which does not contain the alternating group, and every odd prime $p$ such that $\bar X$ is $p$-nilpotent and $p$ divides $|\bar X|$, the program searches deterministically, by increasing size up to $n/2$, for a subset $A$ of $\{1,\dots,n\}$ whose setwise stabilizer has order prime to $p$. Since a subset and its complement have the same stabilizer, this search is exhaustive. There are $49$ such pairs, listed in Table \ref{primtable}. A subset was found in all cases except for $\bar X=\mathrm{A\Gamma L}_1(8)$, \texttt{PrimitiveGroup}$(8,2)$, and $p=3$; for this group the program checks all $3^8$ colourings with three colours, and $4788$ of them are not preserved by any of the $56$ elements of order $3$.

\begin{table}[h]
\footnotesize
\begin{tabular}{rrrrlr|rrrrlr}
$n$ & \# & $|\bar X|$ & $p$ & $A$ & $|\mathrm{St}|$ & $n$ & \# & $|\bar X|$ & $p$ & $A$ & $|\mathrm{St}|$\\ \hline
5 & 1 & 5 & 5 & $\{1\}$ & 1 & 23 & 4 & 506 & 11 & $\{1, 2\}$ & 2\\
7 & 1 & 7 & 7 & $\{1\}$ & 1 & 25 & 1 & 75 & 3 & $\{1, 2\}$ & 1\\
7 & 3 & 21 & 3 & $\{1, 2\}$ & 1 & 25 & 3 & 150 & 3 & $\{1, 2\}$ & 2\\
7 & 4 & 42 & 3 & $\{1, 2\}$ & 2 & 25 & 9 & 300 & 3 & $\{1, 2\}$ & 2\\
8 & 1 & 56 & 7 & $\{1, 2\}$ & 2 & 25 & 14 & 600 & 3 & $\{1, 2\}$ & 2\\
8 & 2 & 168 & 3 & none & --- & 25 & 15 & 600 & 3 & $\{1, 2\}$ & 2\\
11 & 1 & 11 & 11 & $\{1\}$ & 1 & 25 & 18 & 1200 & 3 & $\{1, 2\}$ & 4\\
11 & 3 & 55 & 5 & $\{1, 2\}$ & 1 & 27 & 2 & 351 & 13 & $\{1, 2\}$ & 1\\
11 & 4 & 110 & 5 & $\{1, 2\}$ & 2 & 27 & 6 & 702 & 13 & $\{1, 2\}$ & 2\\
13 & 1 & 13 & 13 & $\{1\}$ & 1 & 29 & 1 & 29 & 29 & $\{1\}$ & 1\\
13 & 3 & 39 & 3 & $\{1, 2\}$ & 1 & 29 & 4 & 203 & 7 & $\{1, 2\}$ & 1\\
13 & 5 & 78 & 3 & $\{1, 2\}$ & 2 & 29 & 5 & 406 & 7 & $\{1, 2\}$ & 2\\
13 & 6 & 156 & 3 & $\{1, 2\}$ & 2 & 29 & 6 & 812 & 7 & $\{1, 2\}$ & 2\\
16 & 1 & 80 & 5 & $\{1, 2\}$ & 2 & 31 & 1 & 31 & 31 & $\{1\}$ & 1\\
16 & 3 & 240 & 3 & $\{1, 2\}$ & 2 & 31 & 3 & 93 & 3 & $\{1, 2\}$ & 1\\
16 & 3 & 240 & 5 & $\{1, 2\}$ & 2 & 31 & 4 & 155 & 5 & $\{1, 2\}$ & 1\\
16 & 6 & 480 & 3 & $\{1, 2\}$ & 4 & 31 & 5 & 186 & 3 & $\{1, 2\}$ & 2\\
17 & 1 & 17 & 17 & $\{1\}$ & 1 & 31 & 6 & 310 & 5 & $\{1, 2\}$ & 2\\
19 & 1 & 19 & 19 & $\{1\}$ & 1 & 31 & 7 & 465 & 3 & $\{1, 2\}$ & 1\\
19 & 3 & 57 & 3 & $\{1, 2\}$ & 1 & 31 & 7 & 465 & 5 & $\{1, 2\}$ & 1\\
19 & 4 & 114 & 3 & $\{1, 2\}$ & 2 & 31 & 8 & 930 & 3 & $\{1, 2\}$ & 2\\
19 & 5 & 171 & 3 & $\{1, 2\}$ & 1 & 31 & 8 & 930 & 5 & $\{1, 2\}$ & 2\\
19 & 6 & 342 & 3 & $\{1, 2\}$ & 2 & 32 & 1 & 992 & 31 & $\{1, 2\}$ & 2\\
23 & 1 & 23 & 23 & $\{1\}$ & 1 & 32 & 2 & 4960 & 5 & $\{1, 2, 3\}$ & 1\\
23 & 3 & 253 & 11 & $\{1, 2\}$ & 1 & & & & & &\\
\end{tabular}
\medskip
\caption{The $49$ pairs $(\bar X,p)$ of Appendix A.1. Here $\bar X$ is \texttt{PrimitiveGroup}$(n,\#)$, $A$ is the first subset (in the order used by \textsf{GAP}) whose setwise stabilizer $\mathrm{St}$ has order prime to $p$, and ``none'' means that no such subset exists.}\label{primtable}
\end{table}

\begin{Verbatim}[fontsize=\footnotesize]
# Appendix A.1: primitive p-nilpotent groups of degree 5..32 (p odd), not containing Alt(n).
# For each pair (G,p) an explicit subset A with p not dividing |Stab_G(A)| is printed,
# found by a deterministic search over subsets ordered by size; "none" means that no such subset exists,
# in which case all 3-colourings are checked.
LoadPackage("primgrp");;
cnt := 0;;
for n in [5..32] do
  for i in [1..NrPrimitiveGroups(n)] do
    G := PrimitiveGroup(n, i);
    if IsNaturalAlternatingGroup(G) or IsNaturalSymmetricGroup(G) then continue; fi;
    for p in Filtered(PrimeDivisors(Size(G)), q -> q > 2) do
      if not IsPNilpotent(G, p) then continue; fi;
      cnt := cnt + 1;
      found := fail;
      for k in [1..QuoInt(n,2)] do
        for A in IteratorOfCombinations([1..n], k) do
          if Size(Stabilizer(G, A, OnSets)) mod p <> 0 then found := A; break; fi;
        od;
        if found <> fail then break; fi;
      od;
      if found <> fail then
        Print(n, " & ", i, " & ", Size(G), " & ", p, " & ", found, " & ", Size(Stabilizer(G, found, OnSets)), "\n");
      else
        el := Filtered(Elements(G), g -> Order(g) = p);
        good := Number(Tuples([1,2,3], n), c -> ForAll(el, g -> Permuted(c, g) <> c));
        Print(n, " & ", i, " & ", Size(G), " & ", p, " & none; good 3-colourings: ", good, "\n");
      fi;
    od;
  od;
od;
Print("pairs: ", cnt, "\n");
QUIT;
\end{Verbatim}

\textbf{A.2. $5$-nilpotent subgroups of $\Sp(6,2)$.} For the proof of Lemma \ref{layers} we determined, using the table of marks of $\Sp(6,2)$, all conjugacy classes of $5$-nilpotent subgroups of $\Sp(6,2)$ of order divisible by $5$. There are six such classes, of orders $5$, $10$, $15$, $30$, $80$ and $160$, and if $H$ is the normal $5$-complement of such a subgroup and $y$ an element of order $5$, then $|H:C_H(y)|$ is $1$ in the first four cases and $16$ in the last two.

\begin{Verbatim}[fontsize=\footnotesize]
# Appendix A.2: 5-nilpotent subgroups of Sp(6,2) of order divisible by 5 (table of marks)
LoadPackage("tomlib");;
tom := TableOfMarks("S6(2)");;
G := UnderlyingGroup(tom);; Print(Size(G), "\n");
ords := OrdersTom(tom);;
best := 0;; cnt := 0;; rec5 := [];;
for i in [1..Length(ords)] do
  if ords[i] mod 5 <> 0 then continue; fi;
  XX := RepresentativeTom(tom, i);
  if not IsPNilpotent(XX, 5) then continue; fi;
  cnt := cnt + 1;
  HH := First(NormalSubgroups(XX), n -> Size(n) = Size(XX)/5);
  y := First(Elements(SylowSubgroup(XX, 5)), g -> Order(g) = 5);
  ind := Index(HH, Centralizer(HH, y));
  Add(rec5, [Size(XX), ind]);
  best := Maximum(best, ind);
od;
Print("5-nilpotent classes: ", cnt, "\n", rec5, "\nmax index: ", best, "\n");
QUIT;
\end{Verbatim}

\textbf{A.3. The group $2.\Sz(8)$.} For Lemma \ref{data3}(c) we determined the degrees of the faithful irreducible ordinary and Brauer characters of $2.\Sz(8)$: they are $40,56,64,104$ in characteristic $0$, $8,40,48,56$ in characteristic $5$, $40,56,64$ in characteristic $7$, and $16,24,40,104$ in characteristic $13$. Moreover, the character table library contains no group $2.\Sz(8).3$, and the centre of $2^2.\Sz(8).3$ is trivial.

\begin{Verbatim}[fontsize=\footnotesize]
# Appendix A.3: 2.Sz(8) -- faithful degrees in characteristics 0, 5, 7, 13, and the outer automorphism
LoadPackage("ctbllib");;
t := CharacterTable("2.Sz(8)");;
z := ClassPositionsOfCentre(t)[2];;
Print("char 0: ", Set(List(Filtered(Irr(t), x -> x[z] <> x[1]), x -> x[1])), "\n");
for p in [5, 7, 13] do
  b := t mod p;;
  zb := Position(GetFusionMap(b, t), z);;
  Print("char ", p, ": ", Set(List(Filtered(Irr(b), x -> x[zb] <> x[1]), x -> x[1])), "\n");
od;
Print("2.Sz(8).3 in the library: ", CharacterTable("2.Sz(8).3") <> fail, "\n");
Print("centre of 2^2.Sz(8).3: ", ClassPositionsOfCentre(CharacterTable("2^2.Sz(8).3")), "\n");
QUIT;
\end{Verbatim}

\textbf{A.4. Centralizers.} The centralizer orders of elements of order $p$ used in the proof of Lemma \ref{layers} were read off from the character table library.

\begin{Verbatim}[fontsize=\footnotesize]
# Appendix A.4: centralizer orders of elements of order p used in Lemma 9.4
LoadPackage("ctbllib");;
for pr in [["S4(2)",5],["S6(2)",5],["S8(2)",5],["S10(2)",5],["S4(3)",5],["S6(2)",7],["S8(2)",7]] do
  t := CharacterTable(pr[1]);
  if t = fail then Print(pr[1], " no table\n"); continue; fi;
  o := OrdersClassRepresentatives(t);; c := SizesCentralizers(t);;
  Print(pr[1], " p=", pr[2], " centralizers of classes of order p: ", List(Filtered([1..Length(o)], i -> o[i]=pr[2]), i -> c[i]), "\n");
od;
QUIT;
\end{Verbatim}

\section{Numerical verifications}

The estimates of Sections 9 and 10 were verified with the following \textsf{Python} programs (Python 3 with \textsf{SymPy}). Each program terminates with an assertion error if one of the inequalities fails; the outputs are included in the supplementary material, and the relevant numbers are reproduced in Tables \ref{exclayers}, \ref{crosschar} and \ref{tight3}.

\textbf{B.1.} \texttt{p5\_layers.py}: the layer inequality in the proof of Lemma \ref{layers} for all primes $p,\ell<200$, which yields exactly the six exceptional layers, and the bounds of Table \ref{exclayers}, including all configurations in which exceptional layers occur together with other active factors.

\begin{Verbatim}[fontsize=\footnotesize]
# Lemma 9.4 (extraspecial layers, p >= 5) and Table 2.
from math import log, floor
from sympy import isprime, factorint, primerange

def sp_order(n, l):
    s = l**(n*n)
    for i in range(1, n+1): s *= l**(2*i)-1
    return s
def minF(p, l):
    """smallest prime power = 1 mod p*l, of characteristic different from l and p (odd if l = 2)."""
    m = p*l; t = 1
    while True:
        x = 1 + m*t; t += 1
        f = factorint(x)
        if len(f) == 1 and list(f)[0] not in (l, p): return x

# (1) the analytic criterion (9.4): (2n^2+3n) log_F(l) <= (1-2/p) l^n - 1, for active layers
exc = []
for p in primerange(5, 200):
    for l in primerange(2, 200):
        if l == p: continue
        F = minF(p, l)
        n = 1
        while l**n < 10**6:
            e = l**n
            active = e >= p-1 and any(pow(l, s, p) in (1, p-1) for s in range(1, n+1))
            if active and (2*n*n+3*n)*log(l)/log(F) > (1-2/p)*e - 1:
                exc.append((p, l, n, e, F))
            n += 1
print("active layers with p, l < 200 violating the layer inequality (proof of Lemma 9.4):", exc)

# (2) exact counts for the six exceptional layers.
# Classes of elements of order p in Sp(2n,l): (number of elements, s) with |[E/Z, y]| = l^(2s);
# centralizer orders from the GAP character table library (see Appendix A.3).
data = {
 (5,2,2): [(4, 2)],                                   # image is C_5 (Lemma 9.3)
 (5,2,3): [(64, 2)],                                  # |H:C_H(y)| <= 16 (Appendix A.2): 16*4 elements
 (5,2,4): [(sp_order(4,2)//3600, 2), (sp_order(4,2)//300, 4)],
 (5,2,5): [(sp_order(5,2)//7257600, 2), (sp_order(5,2)//1800, 4)],
 (5,3,2): [(25920//5, 2)],                            # PSp(4,3) has one class of order 5, centralizer 5
 (7,2,3): [(sp_order(3,2)//7, 3)],
}
print("\n(p,l,e)   |F|   lambda_c  e-1   alone   with another active factor (D>=2)")
res = {}
for (p, l, n), cls in data.items():
    e = l**n; F = minF(p, l)
    N = sum(c * l**(2*s) for c, s in cls)
    lam = log(1+N)/log(F)
    alone = p/(p-1) * sum(c * l**(2*s) * F**(floor((e + (p-1)*l**(n-s))/p) - e) for c, s in cls)
    def T(D, m):   # term of an exceptional factor when m factors are active and D >= 2^(m-1)
        return p/(p-1) * sum(c * l**(2*s) * F**(D*(floor((e+(p-1)*l**(n-s))/p) - e) + D - (m-1)) for c, s in cls)
    res[(p,l,n)] = (F, T)
    print("(%d,%d,%d) %5d  %7.2f  %4d  %.3g   %.3g" % (p, l, e, F, lam, e-1, alone, T(2, 2)))
    assert lam <= e-1 and alone < 0.1

# (3) combinations: m active factors, of which the exceptional ones are listed; the others satisfy (PC).
worst = 0
for p, excs in [(5, [(5,2,2),(5,2,3),(5,2,4),(5,2,5)]), (5, [(5,3,2)]), (7, [(7,2,3)])]:
    for x in excs:
        for with3 in ([None, (5,3,2)] if p == 5 and x[1] == 2 else [None]):
            ex = [x] + ([with3] if with3 else [])
            for m in range(len(ex), 12):
                if m == 1: continue
                D = 2**(m-1)
                F = min(res[y][0] for y in ex)
                tot = sum(res[y][1](max(D, 4 if len(ex) == 2 else D), m) for y in ex)
                tot += (m - len(ex)) * p/(p-1) * F**(1 - m - (D - 1))
                worst = max(worst, tot)
                assert tot < 1, (x, with3, m, tot)
print("\nlargest bound in configurations with an exceptional layer and m >= 2 active factors: %.3g" % worst)
\end{Verbatim}

\textbf{B.2.} \texttt{p5\_blocks.py}: the elementary inequalities used in Lemmas \ref{permfactor} and \ref{fieldfactor}.

\begin{Verbatim}[fontsize=\footnotesize]
# Lemmas 9.5 and 9.7 (Steps 1, 2): elementary inequalities, checked on a range (monotone beyond it).
from math import log, factorial
L8 = lambda x: log(x)/log(8)
# Lemma 9.5, d = 2: 1.6k + log_8 k! <= 0.7*2^k - 1 for k >= 5
assert all(1.6*k + L8(factorial(k)) <= 0.7*2**k - 1 for k in range(5, 200))
# Lemma 9.5, d >= 3: 2k d^2 + log_8 k! <= 0.7 d^k - 1 for k >= 5
assert all(2*k*d*d + L8(factorial(k)) <= 0.7*d**k - 1 for k in range(5, 60) for d in range(3, 60))
# Lemma 9.7, Step 2 (k = number of tensor factors, k >= p >= 5):
#   d_M = 2: 3.5k <= 0.7*2^k - 1;  d_M >= 3: k(d_M^2+2) <= 0.7*d_M^k - 1;  Suzuki/Ree: 54k <= 0.7*4^k - 1
assert all(3.5*k <= 0.7*2**k - 1 for k in range(5, 300))
assert all(k*(d*d+2) <= 0.7*d**k - 1 for k in range(5, 100) for d in range(3, 100))
assert all(54*k <= 0.7*4**k - 1 for k in range(5, 300))
# Lemma 9.7, Step 1: log_8 k! <= 2(2^(k-1) - k) + k - 1
assert all(L8(factorial(k)) <= 2*(2**(k-1) - k) + k - 1 for k in range(1, 300))
print("Lemmas 9.4 and 9.5(a): all inequalities verified")
\end{Verbatim}

\textbf{B.3.} \texttt{p5\_crosschar.py}: the inequality \eqref{ineq} for all classical groups with $n\le8$ and all exceptional groups, at the smallest admissible even and odd $q$ (Table \ref{crosschar}); the case $n\ge9$; and the direct check for $\PSL_2(q)$ and $\Sz(q)$. The bounds $\delta(S)$ used are those of \cite{LS}: $q^{n-1}-1$ for $\PSL_n(q)$; $q(q^{n-1}-1)/(q+1)$ and $(q^n-1)/(q+1)$ for $\PSU_n(q)$ with $n$ odd and even; $(q^m-1)/2$ and $q^{m-1}(q^{m-1}-1)(q-1)/2$ for $\PSp_{2m}(q)$ with $q$ odd and even; $q^{2(m-1)}-1$ for $\Omega_{2m+1}(q)$; $(q^{m-1}\mp1)(q^{m-2}\pm1)$ for $P\Omega^\pm_{2m}(q)$; $q(q^2-1)$, $q^3(q^2-1)$, $q^9(q^2-1)$, $q^{15}(q^2-1)$, $q^{27}(q^2-1)$ for $G_2$, ${}^3D_4$, $E_6$, $E_7$, $E_8$; $q^4(q^6-1)$ and $\frac12q^7(q^3-1)(q-1)$ for $F_4$ with $q$ odd and even; $q^8(q^4+1)(q^3-1)$ for ${}^2E_6$; $q(q-1)$ for ${}^2G_2$; and $(q/2)^{1/2}q^4(q-1)$ for ${}^2F_4$.

\begin{Verbatim}[fontsize=\footnotesize]
# Lemma 9.7, Step 3 (field automorphisms of order p >= 5, cross characteristic) and Table 3.
# Checks the inequality of Step 3 of the proof of Lemma 9.7: (dim G + 2) log_8 q <= 0.8 delta(S)/t - 1, and reports the largest t
# for which it still holds (robustness with respect to the generation bound t).
from math import log, sqrt
L8 = lambda x: log(x)/log(8)
def ls(typ, n, q):
    """Landazuri-Seitz lower bound for the cross-characteristic degree; n = natural dimension / rank."""
    odd = q % 2
    if typ == 'L':  return q**(n-1) - 1
    if typ == 'U':  return q*(q**(n-1)-1)/(q+1) if n % 2 else (q**n-1)/(q+1)
    if typ == 'S':  m = n//2; return (q**m-1)/2 if odd else q**(m-1)*(q**(m-1)-1)*(q-1)/2
    if typ == 'O':  m = (n-1)//2; return q**(2*(m-1)) - 1
    if typ == 'O+': m = n//2; return (q**(m-1)-1)*(q**(m-2)+1)
    if typ == 'O-': m = n//2; return (q**(m-1)+1)*(q**(m-2)-1)
    raise ValueError
def dimG(typ, n):
    if typ in 'LU': return n*n-1
    if typ in ('S', 'O'): m = n//2; return m*(2*m+1)
    if typ in ('O+', 'O-'): m = n//2; return m*(2*m-1)
exc = {  # name: (dim G, delta(q) for even q, delta(q) for odd q, admissible parities)
 'G2':  (14,  lambda q: q*(q*q-1), lambda q: q*(q*q-1), (32, 243)),
 '3D4': (28,  lambda q: q**3*(q*q-1), lambda q: q**3*(q*q-1), (32, 243)),
 'F4':  (52,  lambda q: q**7*(q**3-1)*(q-1)/2, lambda q: q**4*(q**6-1), (32, 243)),
 'E6':  (78,  lambda q: q**9*(q*q-1), lambda q: q**9*(q*q-1), (32, 243)),
 '2E6': (78,  lambda q: q**8*(q**4+1)*(q**3-1), lambda q: q**8*(q**4+1)*(q**3-1), (32, 243)),
 'E7':  (133, lambda q: q**15*(q*q-1), lambda q: q**15*(q*q-1), (32, 243)),
 'E8':  (248, lambda q: q**27*(q*q-1), lambda q: q**27*(q*q-1), (32, 243)),
 '2G2': (14,  None, lambda q: q*(q-1), (243,)),
 '2F4': (52,  lambda q: q**4*(q-1)*sqrt(q/2), None, (32,)),
}
rows = []
def row(name, q, dG, delta, t):
    lhs = (dG + 2)*L8(q); rhs = 0.8*delta/t - 1
    tmax = 0.8*delta/(lhs + 1)
    rows.append((name, q, lhs, rhs, tmax)); assert lhs <= rhs, (name, q)
for typ, nmin in [('L',3), ('U',3), ('S',4), ('O',7), ('O+',8), ('O-',8)]:
    for n in range(nmin, 9 if typ not in ('S','O+','O-') else 9, 1 if typ in 'LU' else 2):
        for q in (32, 243):
            if typ == 'O' and q % 2 == 0: continue
            row('%s_%d' % (typ, n), q, dimG(typ, n), ls(typ, n, q), n+3)
for name, (dG, de, do, qs) in exc.items():
    for q in qs:
        row(name, q, dG, (de if q % 2 == 0 else do)(q), 11)
print("%-6s %5s %8s %12s %10s" % ("type", "q", "LHS", "RHS", "max t"))
for r in rows: print("%-6s %5d %8.1f %12.4g %10.4g" % r)
print("smallest admissible t over the table: %.1f" % min(r[4] for r in rows))
# classical groups with natural dimension n >= 9: delta >= (q^(n/2-1)-1)/2, dim G <= n^2
for n in range(9, 200):
    for q in (32, 243):
        assert (n*n+2)*L8(q) <= 0.8*(q**(n/2-1)-1)/2/(2*n) - 1
print("n >= 9: inequality holds with t = 2n for 9 <= n < 200 (and monotone beyond)")
# PSL_2(q) and Sz(q): t = 2, kappa >= (1-1/p)/2, checked directly (|F| >= 11 if p = 5, q >= 32; |F| >= 8, q >= 128 if p >= 7)
for p, F, qs in [(5, 11, [32, 243, 1024, 3125]), (7, 8, [128, 2187]), (11, 23, [2048])]:
    k = (1-1/p)/2
    for q in qs:
        f = round(log(q, 2)) if q % 2 == 0 else round(log(q, 3))
        autL2 = q**3*f*2; dL2 = (q-1)/(1 if q % 2 == 0 else 2)
        assert log(autL2)/log(F) <= k*dL2 - 1, ('L2', p, q)
        if q % 2 == 0:
            autSz = q*q*(q*q+1)*(q-1)*f; dSz = (q-1)*sqrt(q/2)
            assert log(autSz)/log(F) <= k*dSz - 1, ('Sz', p, q)
print("PSL_2(q) and Sz(q): (PC) verified at the smallest admissible q for p = 5, 7, 11")
\end{Verbatim}

\textbf{B.4.} \texttt{p3\_sweep.py}: the verification of $\Theta<1$ in the proof of Proposition \ref{p3} for all configurations with $P<2^{10}$, and the elementary inequality at the end of the proof of Lemma \ref{tail3}.

\begin{Verbatim}[fontsize=\footnotesize]
# Verification of Proposition 10.3 (p = 3, K nonsolvable).
# Every configuration of tensor factors with P = prod d_c < 2^10 is checked against the criterion
#   Theta = (3/2) * sum_{c active} N_c * prod_{c' != c, active} (1+N_{c'}) * |F|^(-kappa_c P) < 1,
# with the data of Lemma 10.1.  Configurations with P >= 2^10 are covered by Lemma 10.2.
# The script prints log_|F|(Theta) for the tightest configurations; all values must be negative.
from math import log, prod, factorial, sqrt, ceil
from sympy import factorint, isprime

def lg(x, b): return log(x)/log(b)
def sp_order(n, l):
    s = l**(n*n)
    for i in range(1, n+1): s *= l**(2*i)-1
    return s
def aut_sz(q):              # |Aut(Sz(q))| = q^2 (q^2+1)(q-1) f
    f = round(log(q, 2)); return q*q*(q*q+1)*(q-1)*f
LIM = 2**10

# ---------- factor types ----------
# each factor: dict(d, N, kappa, active, kind, data)
def layer(l, n, sz=False):
    d = l**n
    if l == 2: kap = 0.5
    else:      kap = 1 - (1 + 2/l)/3
    if sz and l == 2 and n == 6:
        N = 2**12 * 3 * 7 * 29120            # Lemma 10.1(f): normalizer of Sz(8) in Sp(12,2)
    else:
        N = l**(2*n) * sp_order(n, l)        # all of l^{2n}.Sp(2n,l)
    return dict(d=d, N=N, kappa=kap, active=True, kind=('L', l, n, sz))
layers = []
for l in [x for x in range(2, 257) if isprime(x) and x != 3]:
    n = 1
    while l**n < LIM:
        layers.append(layer(l, n))
        n += 1
# Suzuki-carrying layers (case E(K)=1): l=2, n>=6; odd l with l^n >= 121
szlayers = [layer(2, 6, sz=True)] + [layer(2, n) for n in range(7, 10)] + \
           [layer(l, n) for l in [5,7,11,13,17,19,23,29,31] for n in range(1,6) if 121 <= l**n < LIM]

def szblock(q, js, char2, cover=False, kcross=1, char5=False):
    """A block of components Sz(q) (or 2.Sz(8)).  char2: the modules have dimensions 4^j, j in js
    (Steinberg tensor products, Lemma 10.1(d)); otherwise kcross components in cross characteristic,
    each of the minimal dimension allowed by Landazuri-Seitz (larger dimensions only help)."""
    f = round(log(q, 2))
    szo = lambda x: x*x*(x*x+1)*(x-1)
    if char2:
        k = len(js); dims = [4**j for j in js]
        kfa = None
        for j in js:
            if j % 3 == 0 and f % 3 == 0:               # linear field automorphism: alpha-invariant module
                d1 = 4**(j//3); x = 1 - (d1**3 + 2*d1)/(3*d1**3)
                kfa = x if kfa is None else min(kfa, x)
        kperm = None
        for j in set(js):
            if js.count(j) >= 3:
                d = 4**j; x = 1 - (d**3 + 2*d)/(3*d**3)
                kperm = x if kperm is None else min(kperm, x)
    elif cover:
        # 2.Sz(8) (Lemma 10.1(c)): no automorphism of order 3 (the outer automorphism of Sz(8) permutes
        # the three double covers); faithful modules have dimension >= 8 if r = 5 and >= 16 otherwise.
        k = kcross; d = 8 if char5 else 16; dims = [d]*k
        kfa = None
        kperm = 1 - (d**3 + 2*d)/(3*d**3) if k >= 3 else None
    else:
        k = kcross; d = ceil((q-1)*sqrt(q/2)); dims = [d]*k
        kfa = 1/3 if f % 3 == 0 else None
        kperm = 1 - (d**3 + 2*d)/(3*d**3) if k >= 3 else None
    opts = [x for x in (kfa, kperm) if x is not None]
    if char2:
        # |Aut Sz(q)| < q^6 <= |F|^(6 j) for each component; k = 1: 2|Sz(q):Sz(q^(1/3))| < 2.2 q^(10/3)
        jm = max(js)
        if k == 1:
            Nf = (lambda F, jm=jm: 2.2 * F**(10*jm/3)) if kfa is not None else (lambda F: 0)
        else:
            Nf = lambda F, js=tuple(js), k=k: prod(F**(6*j) for j in js) * factorial(k)
    else:
        A = szo(8) if cover else aut_sz(q)        # Aut(2.Sz(8)) = Inn(2.Sz(8)) has order |Sz(8)|
        N = A**k * factorial(k)
        if k == 1:
            N = 2 * szo(q) // szo(round(q**(1/3))) if kfa is not None else 0
        Nf = lambda F, N=N: N
    N = 1
    return dict(d=prod(dims), N=N, Nf=Nf, kappa=min(opts) if opts else None, active=bool(opts),
                kind=('Sz', q, tuple(js) if char2 else kcross, char2, cover, char5), q=q, js=js, char2=char2, char5=char5)

szblocks = []
def jtuples(maxsum):
    out = []
    def rec(cur, mn, sm):
        if cur: out.append(list(cur))
        for j in range(mn, maxsum - sm + 1):
            rec(cur + [j], j, sm + j)
    rec([], 1, 0)
    return out
for js in jtuples(4):
    szblocks.append(szblock(2**3 if max(js) <= 3 else 2**max(js) if max(js) % 2 else 2**(max(js)+1), js, True))
for q in [8, 32, 128]:
    for k in range(1, 3):
        b = szblock(q, None, False, kcross=k)
        if b['d'] < LIM: szblocks.append(b)
for c5 in (True, False):
    for k in range(1, 4):
        b = szblock(8, None, False, cover=True, kcross=k, char5=c5)
        if b['d'] < LIM: szblocks.append(b)

def prime_powers(maxv=10**7):
    for x in range(2, 200000):
        fct = factorint(x)
        if len(fct) == 1: yield x, list(fct)[0]


def min_field(conf):
    """smallest |F| compatible with the configuration (Lemma 10.1(a),(d))."""
    mod = 3
    forbid = {3}
    char = None; need = 1
    for c in conf:
        if c['kind'][0] == 'L':
            l = c['kind'][1]; mod *= l; forbid.add(l)
            if l == 2: char = 'odd' if char in (None, 'odd') else ('X' if char != 5 else 5)
        else:
            if c['char2']:
                char = 2 if char in (None, 2) else 'X'
                if 1 in c['js']: need = 6 if need % 6 else need          # GF(q) <= F with q >= 8, and 3 | |F|-1
            elif c.get('char5'):
                char = 5 if char in (None, 'odd', 5) else 'X'
            else:
                char = 'odd' if char in (None, 'odd') else ('X' if char != 5 else 5)
    if char == 'X': return None
    if char == 5:
        if 5 in forbid: return None
        g = 1
        while (5**g - 1) % mod: g += 1
        return 5**g
    if char == 2:
        g = need
        while (2**g - 1) % mod: g += need
        return 2**g
    t = 1
    while True:
        x = 1 + mod*t; t += 1
        fct = factorint(x)
        if len(fct) != 1: continue
        r = list(fct)[0]
        if r in forbid: continue
        if char == 'odd' and r == 2: continue
        if char == 5 and r != 5: continue
        return x


def check(conf):
    F = min_field(conf)
    if F is None: return None
    P = prod(c['d'] for c in conf)
    act = [c for c in conf if c['active']]
    if not act: return -1e9                           # Corollary 7.3: nothing to check
    for c in act:
        if 'Nf' in c: c['N'] = c['Nf'](F)
    lam = [lg(1 + c['N'], F) for c in act]  # base |F| logarithms
    L = sum(lam)
    ex = [L - lam[i] + lg(max(c['N'], 1), F) - c['kappa']*P for i, c in enumerate(act)]
    # monotonicity in |F|: the power of |F| in every term must be negative (proof of Prop. 10.3)
    Fpow = sum(6*sum(c['js']) if len(c['js']) > 1 else 10*max(c['js'])/3
               for c in act if c['kind'][0] == 'Sz' and c['char2'])
    for c in act: assert Fpow - c['kappa']*P < 0, ([x['kind'] for x in conf], Fpow, c['kappa']*P)
    M = max(ex)
    return M + lg(1.5*sum(__import__("math").exp((e-M)*log(F)) for e in ex), F)  # log_|F| of the bound; must be < 0

worst = None; count = 0; fails = []; allv = []
carriers = szlayers + szblocks
others_layers = layers
# configurations: one or more Suzuki-carrying factors, plus layers with distinct primes
def layer_sets(maxprod):
    out = [[]]
    by_l = {}
    for c in others_layers: by_l.setdefault(c['kind'][1], []).append(c)
    ls = list(by_l)
    def rec(i, cur, pr):
        if i == len(ls):
            out.append(list(cur)); return
        rec(i+1, cur, pr)
        for c in by_l[ls[i]]:
            if pr * c['d'] < maxprod: rec(i+1, cur+[c], pr*c['d'])
    out = []
    rec(0, [], 1)
    return out
LS = layer_sets(LIM)
carriers.sort(key=lambda c: c['d'])
def keyof(c): return ('L', c['kind'][1]) if c['kind'][0] == 'L' else ('Sz', c['q'])
def rec(start, cur, pr, keys):
    global count, worst
    if cur:
        for ls in LS:
            if any(('L', x['kind'][1]) in keys for x in ls): continue
            if pr * prod(x['d'] for x in ls) >= LIM: continue
            conf = cur + ls
            v = check(conf)
            if v is None: continue
            count += 1
            allv.append((v, [c['kind'] for c in conf], min_field(conf), prod(c['d'] for c in conf)))
            if worst is None or v > worst[0]: worst = (v, [c['kind'] for c in conf])
            if v >= 0: fails.append((v, [c['kind'] for c in conf]))
    for i in range(start, len(carriers)):
        c = carriers[i]
        if pr * c['d'] >= LIM: break
        if c['kind'][0] == 'Sz' and c['char2']:
            rec(i, cur + [c], pr * c['d'], keys)        # several blocks Sz(q), Sz(q'), q != q', may share j's
            continue
        if keyof(c) in keys: continue
        rec(i+1, cur + [c], pr * c['d'], keys | {keyof(c)})
rec(0, [], 1, frozenset())
print("configurations checked:", count)
print("failures:", len(fails))
for x in fails[:20]: print(x)
print("largest value of log_|F|(bound):", worst)
allv.sort(key=lambda x: -x[0])
print("ten tightest configurations (value, factors, |F|, P):")
for x in allv[:10]: print("  %.2f" % x[0], x[1], x[2], x[3])
# tail: Lemma 10.2
for L in range(10, 41):
    lhs = 2*L*L + 5*L + L*log(L,2) + log(1.5*L, 4)
    assert 2**L/3 > lhs, L
print("tail inequality 2^L/3 > 2L^2+5L+L log2 L + log_4(1.5L) verified for 10 <= L <= 40 (and it is monotone beyond)")

# ---- export of all checked configurations (canonical form) for the independent check ----
import json
def canon(kinds):
    out = []
    for k in kinds:
        if k[0] == 'L': out.append(('L', k[1], k[2], bool(k[3])))
        elif k[3]:      out.append(('C2', tuple(sorted(k[2]))))
        else:           out.append(('X', k[1], k[2], bool(k[4]), bool(k[5])))
    return tuple(sorted(out, key=repr))
exp = {}
for v, kinds, F, P in allv:
    exp[repr(canon(kinds))] = v
json.dump(exp, open('p3_sweep_configs.json', 'w'))
print("distinct configurations:", len(exp))
\end{Verbatim}

\textbf{B.5.} \texttt{p3\_independent.py}: an independent implementation of the computation of B.4, which compares its results with those of B.4. For the monotonicity test, it computes for each configuration the three successive admissible values of $|\F|$ following the smallest one (subject to the same conditions as the smallest one), and compares $\Theta$ at each of them with the value at the preceding one; it also recomputes $\Theta$ with larger dimensions of the components in odd characteristic. Its output is: $713$ canonical configurations, of which $386$ contain active factors and a compatible field; largest value of $\log_{|\F|}\Theta$ equal to $-0.35914$; no disagreement with B.4; no configuration missing in either program; $1158=3\cdot386$ field sizes tested, and no violation of monotonicity. For example, for a single component $\Sz(8)$ in odd characteristic the fields tested are $\GF(7)$, $\GF(13)$, $\GF(19)$ and $\GF(25)$, with $\log_{|\F|}\Theta\le-0.359$, $-1.399$, $-1.820$ and $-2.063$, respectively.

\begin{Verbatim}[fontsize=\footnotesize]
# Independent re-implementation of the computation in the proof of Proposition 10.3 (p = 3).
# Written separately from p3_sweep.py: different enumeration (canonical multisets of factor types),
# high-precision arithmetic (mpmath), an independent search for the smallest admissible |F|, and
# additional tests of the monotonicity claims.  It compares its results with p3_sweep_configs.json.
import json
from mpmath import mp, mpf, log as mlog, sqrt as msqrt, ceil as mceil
from sympy import primerange, factorint
mp.dps = 50
LIM = 1024

def sp(n, l):
    r = l**(n*n)
    for i in range(1, n+1): r *= l**(2*i) - 1
    return r
def szo(q): return q*q*(q*q+1)*(q-1)

# ---------------- factor types (Lemma 10.1) ----------------
# each factor: (key, d, kappa or None, Nfun(F), carrier, constraints)
F = {}
def L(l, n, sz=False):
    d = l**n
    kap = mpf(1)/2 if l == 2 else 1 - (1 + mpf(2)/l)/3
    N = 2**12*3*7*szo(8) if sz else l**(2*n)*sp(n, l)
    carrier = sz or (l == 2 and n >= 7) or (l > 2 and n >= 2 and d >= 121)   # Lemma 10.1(f)
    return dict(key=('L', l, n, sz), d=d, kap=kap, N=lambda Fq, N=N: mpf(N), carrier=carrier,
                ell=l, char=('odd' if l == 2 else None))
def X(q, k, cover, c5=False):
    f = {8: 3, 32: 5, 128: 7}[q]
    from math import factorial
    if cover:
        # 2.Sz(8): Aut = Inn (order |Sz(8)|), no element of order 3 acts on a single component;
        # faithful modules: dimension >= 8 in characteristic 5, >= 16 otherwise (Brauer tables)
        d1 = 8 if c5 else 16
        kaps = [1 - mpf(d1**3 + 2*d1)/(3*d1**3)] if k >= 3 else []
        N = szo(8)**k*factorial(k) if k >= 3 else 0
        return dict(key=('X', q, k, True, c5), d=d1**k, kap=min(kaps) if kaps else None,
                    N=lambda Fq, N=N: mpf(N), carrier=True, ell=None, char=(5 if c5 else 'odd'), qkey=q)
    d1 = int(mceil((q-1)*msqrt(mpf(q)/2)))
    kaps = []
    if f % 3 == 0: kaps.append(mpf(1)/3)
    if k >= 3: kaps.append(1 - mpf(d1**3 + 2*d1)/(3*d1**3))
    aut = szo(q)*f*(2 if cover else 1)
    from math import factorial
    if k == 1: N = 2*szo(q)//szo(round(q**(1/3))) if f % 3 == 0 else 0
    else:      N = aut**k*factorial(k)
    return dict(key=('X', q, k, False, False), d=d1**k, kap=min(kaps) if kaps else None,
                N=lambda Fq, N=N: mpf(N), carrier=True, ell=None, char='odd', qkey=q)
def C2(js):
    from math import factorial
    js = tuple(sorted(js)); k = len(js)
    kaps = []
    for j in js:
        if j % 3 == 0:
            d1 = 4**(j//3); kaps.append(1 - mpf(d1**3 + 2*d1)/(3*d1**3))
    for j in set(js):
        if js.count(j) >= 3:
            d = 4**j; kaps.append(1 - mpf(d**3 + 2*d)/(3*d**3))
    if k == 1:
        Nf = (lambda Fq, j=js[0]: mpf('2.2')*mpf(Fq)**(mpf(10*j)/3)) if kaps else (lambda Fq: mpf(0))
    else:
        Nf = lambda Fq, js=js, k=k: mpf(Fq)**(6*sum(js))*factorial(k)
    d = 1
    for j in js: d *= 4**j
    return dict(key=('C2', js), d=d, kap=min(kaps) if kaps else None, N=Nf, carrier=True,
                ell=None, char=2, need6=(1 in js))

layers = []
for l in primerange(2, 1024):
    if l == 3: continue
    n = 1
    while l**n < LIM:
        layers.append(L(l, n)); n += 1
layers.append(L(2, 6, sz=True))
xs = [X(q, k, False) for q in (8, 32, 128) for k in (1, 2)
      if int(mceil((q-1)*msqrt(mpf(q)/2)))**k < LIM]
xs += [X(8, k, True, c5) for c5 in (True, False) for k in (1, 2, 3) if (8 if c5 else 16)**k < LIM]
def parts(maxsum):
    res = []
    def r(cur, mn, s):
        if cur: res.append(tuple(cur))
        for j in range(mn, maxsum - s + 1): r(cur + [j], j, s + j)
    r([], 1, 0); return res
c2s = [C2(js) for js in parts(4)]

# ---------------- smallest admissible field ----------------
def is_pp(x):
    f = factorint(x); return (list(f)[0] if len(f) == 1 else None)
def field(conf):
    chars = {c['char'] for c in conf if c['char'] is not None}
    if chars == {'odd', 5}: chars = {5}
    if len(chars) > 1: return None
    ch = chars.pop() if chars else None
    M = 3; forb = {3}
    for c in conf:
        if c['ell']: M *= c['ell']; forb.add(c['ell'])
    need6 = any(c.get('need6') for c in conf)
    if ch == 2:
        g = 1
        while True:
            if (2**g - 1) % M == 0 and (not need6 or g % 6 == 0): return 2**g
            g += 1
    if ch == 5:
        if 5 in forb: return None
        g = 1
        while (5**g - 1) % M: g += 1
        return 5**g
    x = 1
    while True:
        x += M
        p = is_pp(x)
        if p is None or p in forb: continue
        if ch == 'odd' and p == 2: continue
        if ch == 5 and p != 5: continue
        return x

def next_field(conf, cur):
    """the smallest admissible |F| (conditions of field()) which is larger than cur."""
    chars = {c['char'] for c in conf if c['char'] is not None}
    if chars == {'odd', 5}: chars = {5}
    ch = chars.pop() if chars else None
    M = 3; forb = {3}
    for c in conf:
        if c['ell']: M *= c['ell']; forb.add(c['ell'])
    need6 = any(c.get('need6') for c in conf)
    if ch == 2:
        g = cur.bit_length()          # cur = 2^(g-1), so start with the next exponent
        while True:
            if (2**g - 1) % M == 0 and (not need6 or g % 6 == 0): return 2**g
            g += 1
    if ch == 5:
        g = 1
        while 5**g <= cur or (5**g - 1) % M: g += 1
        return 5**g
    x = cur
    while True:
        x += M
        p = is_pp(x)
        if p is None or p in forb: continue
        if ch == 'odd' and p == 2: continue
        if ch == 5 and p != 5: continue
        return x

def theta(conf, Fq, dims=None):
    P = 1
    for i, c in enumerate(conf): P *= (dims[i] if dims else c['d'])
    act = [c for c in conf if c['kap'] is not None and c['N'](Fq) > 0]
    if not act: return None
    Ns = [c['N'](Fq) for c in act]
    tot = mpf(0)
    for i, c in enumerate(act):
        t = mpf('1.5')*Ns[i]
        for j in range(len(act)):
            if j != i: t *= 1 + Ns[j]
        tot += t*mpf(Fq)**(-c['kap']*P)
    return tot

# ---------------- enumeration of canonical configurations ----------------
pool = layers + xs + c2s
pool.sort(key=lambda c: (c['d'], repr(c['key'])))
confs = []
def rec(i, cur, prod_):
    if cur and any(c['carrier'] for c in cur): confs.append(list(cur))
    for j in range(i, len(pool)):
        c = pool[j]
        if prod_*c['d'] >= LIM: break
        if c['key'][0] == 'L' and any(x['key'][0] == 'L' and x['ell'] == c['ell'] for x in cur): continue
        if c['key'][0] == 'X' and any(x['key'][0] == 'X' and x['qkey'] == c['qkey'] for x in cur): continue
        rec(j + (0 if c['key'][0] == 'C2' else 1), cur + [c], prod_*c['d'])
rec(0, [], 1)

ref = json.load(open('p3_sweep_configs.json'))
def canon(conf):
    out = []
    for c in conf:
        k = c['key']
        if k[0] == 'L':  out.append(('L', k[1], k[2], bool(k[3])))
        elif k[0] == 'X': out.append(('X', k[1], k[2], bool(k[3]), bool(k[4])))
        else:            out.append(('C2', tuple(sorted(k[1]))))
    return repr(tuple(sorted(out, key=repr)))

checked = 0; worst = (-1e9, None); mism = 0; missing = []; mono_fail = 0; fields_tested = 0; traces = []
for conf in confs:
    Fq = field(conf)
    if Fq is None: continue
    th = theta(conf, Fq)
    key = canon(conf)
    if th is None:
        continue
    checked += 1
    v = float(mlog(th)/mlog(Fq))
    assert v < 0, (key, v)
    if v > worst[0]: worst = (v, key)
    if key in ref:
        if abs(ref[key] - v) > 1e-9: mism += 1
    else:
        missing.append(key)
    # monotonicity tests: (i) the next three admissible fields, compared successively,
    # (ii) larger dimensions in odd characteristic
    cur, prev = Fq, th
    trace = [(Fq, float(mlog(th)/mlog(Fq)))]
    for t in range(3):
        cur = next_field(conf, cur)
        th2 = theta(conf, cur)
        if th2 is not None and th2 > prev*(1 + mpf(10)**-30): mono_fail += 1
        if th2 is not None: prev = th2
        fields_tested += 1
        if th2 is not None: trace.append((cur, float(mlog(th2)/mlog(cur))))
    traces.append((key, trace))
    if any(c['key'][0] == 'X' for c in conf):
        dims = [c['d'] + (5 if c['key'][0] == 'X' else 0) for c in conf]
        th3 = theta(conf, Fq, dims)
        if th3 is not None and th3 > th: mono_fail += 1

print("canonical configurations enumerated:", len(confs))
print("configurations with active factors and a compatible field:", checked)
print("largest log_|F|(Theta): %.6f" % worst[0], worst[1])
print("disagreements with p3_sweep.py (|difference| > 1e-9):", mism)
print("configurations with active factors not present in p3_sweep.py:", len(missing))
print("configurations of p3_sweep.py not found here:", len(set(ref) - {canon(c) for c in confs}))
print("field sizes tested for monotonicity:", fields_tested)
print("monotonicity violations (next three admissible |F|, larger dimensions):", mono_fail)
print("the five tightest configurations, with log_|F|(Theta) at the smallest and the next three fields:")
tr = dict(traces)
for k in sorted(tr, key=lambda k: -tr[k][0][1])[:5]:
    print("  ", k, ["|F|=%d: %.3f" % x for x in tr[k]])
# sanity check of next_field
print("example: admissible fields after 7 for a single Sz(8) block:",
      [7] + [next_field([xs[0]], x) for x in (7, 13, 19)])
\end{Verbatim}

\section*{}
\begin{center}
    \textbf{Acknowledgement}
\end{center}
 Thie author was supported by a grant from the Simons Foundation (SFI-MPS-TSM-00014153 to Thomas M. Keller).

\section*{AI statement}

This paper is a natural sequel to \cite{K}. After \cite{K} was completed, the author used Claude (Opus 5.5) by Anthropic in a lengthy and extended interaction to generalize its results to arbitrary finite groups. The author posed the problem and proposed the overall strategy of the reduction. The AI system developed most of the arguments of Sections 3 to 12, including Theorem \ref{galois} and the examples of Section 12, carried out the numerical estimates and the computer calculations, and produced a first version of the manuscript, which was edited by the author. All arguments and computations were subsequently checked by the author, who takes full responsibility for the content of this paper.

\section*{Disclosure statement}

The author reports no conflict of interest.

\section*{Data availability}

No data were used in this work.

\end{document}